\documentclass[11pt]{article}

\usepackage{amsfonts}

\usepackage{CJK}
\usepackage{amsfonts,amssymb,amsmath,mathrsfs,multirow,booktabs}
\usepackage{xcolor,soul}
\usepackage{color,latexsym,amsfonts}
\newcommand{\qed}{\hfill\rule{2mm}{3mm}\vspace{4mm}}

\newtheorem{theorem}{Theorem}[section]
\newtheorem{lemma}[theorem]{Lemma}
\newtheorem{corollary}[theorem]{Corollary}
\newtheorem{proposition}[theorem]{Proposition}
\newtheorem{example}[theorem]{Example}
\newtheorem{Definition}[theorem]{Definition}
\newtheorem{remark}[theorem]{Remark}
\newtheorem{condition}[theorem]{Condition}
\newtheorem{conjecture}[theorem]{Conjecture}

\def\blemma{\begin{lemma}\sl{}\def\elemma{\end{lemma}}}
\def\btheorem{\begin{theorem}\sl{}\def\etheorem{\end{theorem}}}

\def\bdefinition{\begin{Definition}\sl{}\def\edefinition{\end{Definition}}}
\def\bproposition{\begin{proposition}\sl{}\def\eproposition{\end{proposition}}}
\def\bremark{\begin{remark}\sl{}\def\eremark{\end{remark}}}
\def\bcondition{\begin{condition}\sl{}\def\econdition{\end{condition}}}

\def\beqlb{\begin{eqnarray}}\def\eeqlb{\end{eqnarray}}
\def\beqnn{\begin{eqnarray*}}\def\eeqnn{\end{eqnarray*}}

\def\mbf{\mathbf}

\def\<{\langle}\def\>{\rangle}

\def\ar{&\!\!}

\def\eqref#1{{\rm(\ref{#1})}}

\def\proof{\noindent{\it
		Proof.~}}\def\qed{\hfill$\Box$\medskip}

\def\blue{\color{blue}}

\let\blue\relax

\def\e{{\mbox{\rm e}}}
\def\<{\left<}\def\>{\right>}

\def\mbf{\mathbf}
\newcommand{\dd}{\mathrm{d}}

\newfam\msbmfam\font\tenmsbm=msbm10\textfont
\msbmfam=\tenmsbm\font\sevenmsbm=msbm7
\scriptfont\msbmfam=\sevenmsbm

\def\<{\left<}\def\>{\right>}
\def\({\left(}\def\){\right)}

\begin{document}

	\centerline{\Large\bf Extinction behaviour for mutually}
	\medskip
	\centerline{\Large\bf enhancing continuous-state population dynamics}

	\bigskip
	
	\centerline{
		Jie Xiong\footnote{Department of Mathematics and SUSTech International center for Mathematics, Southern University of Science \& Technology, Shenzhen, China.
			Supported by National Key R\&D Program of China (No.~2022YFA1006102) and
			National Natural Science Foundation of China Grant 12471418.  Email: xiongj@sustech.edu.cn},
		Xu Yang\footnote{School of Mathematics and Information Science,
			North Minzu University, Yinchuan, China. Supported by
			NSFC (No.~12471135). Email: xuyang@mail.bnu.edu.cn. Corresponding author.} and Xiaowen Zhou\footnote{Department of Mathematics and Statistics, Concordia
			University, Montreal, Canada.
			Supported by NSERC (RGPIN-2021-04100).
			Email: xiaowen.zhou@concordia.ca.}}

	\bigskip\bigskip
	
	{\narrower{\narrower
			
			\noindent{\bf Abstract.}
In this paper, we study a two-dimensional
process arising as the unique nonnegative solution to a system of two
stochastic differential equations (SDEs) with mutually enhancing two-way interactions driven by independent Brownian motions
and spectrally positive $\alpha$-stable random measures.
Such a SDE system can be identified as a continuous-state Lotka-Volterra type population model. Extinction properties of the populations are studied for different choices of the coefficients involved in the SDEs.

			\bigskip
			
			
			\bigskip
			
			\par}\par}
	
	\section{Introduction and main results}\label{intro}
	
	\setcounter{equation}{0}

Continuous-state branching processes (CSBPs for short)
are mathematical models that describe the evolution of rescaled populations.
They arise as scaling limits of Galton-Watson branching processes
and find applications in various fields such as biology, population genetics, physics, chemistry, and so on.
We refer to \cite{Ky,Li,Li12} for reviews
and literature on CSBPs.
Extinction behavior is a key topic in the study of population models.
Through the Lamperti random time change, the CSBP
is associated with a spectrally positive L\'evy process stopped when first reaching $0$, which allows many
semiexplicit expressions.
In particular, a sufficient and necessary condition, called Grey's condition, is obtained  in \cite{Grey74}.

In recent years, many authors have studied the generalized nonlinear CSBP characterized
by the following SDE.
\beqlb\label{1.0}
X_t=
X_0+\int_0^t \gamma_0(X_s)\dd s+\int_0^t\int_0^{\gamma_1(X_s)}W(\dd s,\dd u)
+\int_0^t\int_0^\infty\int_0^{\gamma_2(X_{s-})}z\tilde{M}(\dd s,\dd z,\dd u),
 \eeqlb
where $X_0>0$ and $\gamma_0$ and $\gamma_1,\gamma_2\ge0$ are Borel functions on $[0,\infty)$,
$W$ and $\tilde{N}$ are Gaussian white noise and compensated Possion random measure.
If $\gamma_i(x)=\gamma_ix$ for some constant $\gamma_i$, then the solution to \eqref{1.0} reduces to
the classical CSBP.
For $\gamma_2\equiv0$ and $\gamma_1(x)=\gamma_1x$, \eqref{1.0} is also called the
Dawson-Li SDE with its solution  called a CSBP with competition in \cite{Pardoux16},
where the function $\gamma_0$ models an impact of the current population size on the  reproduction dynamics of individuals.

Moreover, if $\gamma_0(x)=c_1x-c_2x^2$ for $c_1,c_2>0$,
then the solution to SDE \eqref{1.0}, called
the branching process with logistic
growth" or ``logistic branching process" in \cite{Lambert2005},  models the density dependence in population.
This density dependence corresponds to intraspecific competition pressure, which is ubiquitous in ecology, and translates mathematically into a quadratic death rate.
The terms $\gamma_i(x)/x$ for $i=1,2$ can be interpreted as population-size-dependent branching rates, and the drift  involving $\gamma_0$ can also be related to continuous population-size-dependent immigration.
Population-size-dependent CSBPs arise as scaling limits of the corresponding
discrete-state branching processes (see, e.g., \cite{LiY06,LiY09}).
If $\gamma_i(x)=\gamma_ix^r$ for $i=0,1,2$, the solution to SDE
\eqref{1.0} is called a CSBP with polynomial branching whose expectation of  extinction time is discussed via Lamperti type transformations in \cite{LiP19},
where the parameter $r$ describes the degree of interaction.
The CSBPs with polynomial branching arise as time-space scaling limits of discrete-state nonlinear branching processes. Intuitively, the functions $\gamma_1$ and $\gamma_2$ are population-dependent rates that generate a small and large number of branching events in children, respectively.
Rather sharp conditions on extinction/non-extinction for SDE \eqref{1.0}
are given in \cite{LYZh} using martingale techniques, and critical cases are further studied in \cite{MYZ21}.
The exponential ergodicity and the strong exponential ergodicity for the
solution to SDE \eqref{1.0} is discussed in \cite{LiWang}
 using coupling techniques.
More recently,  models related to \eqref{1.0} have also
attracted attention.
The extinction/non-extinction for generalized  nonlinear CSBPs with catastrophes
is studied in \cite{Marguet}.
\cite{baiyang} establishes
necessary and sufficient conditions of extinction for the generalized  nonlinear CSBPs with Neveu's branching.
The ergodic property of a CSBP with immigration and competition is given
in \cite{WangLiLiZhou}.
General  nonlinear CSBPs with competition in the L\'{e}vy  environment
and CSBPs in varying environments are studied in \cite{MZ24} and \cite{Fangli}, respectively.

By \cite[Theorems 7.1 and 8.1]{Li12},
on an enlarged probability space,
on the right of \eqref{1.0}, the third term can be transformed into
$\int_0^t\sqrt{\gamma_1(X_s)}\dd B_s$, and the fourth term
can be converted into $\int_0^t\int_0^\infty z\tilde{M}_{\gamma_2}(\dd s,\dd z)$,
where $(B_t)_{t\ge0}$ is a Brownian motion and $\tilde{M}_{\gamma_2}(\dd s,\dd z)$
is an optional compensated Poisson random measure.
Using the Lamperti transform for positive self-similar Markov processes,
in certain cases, \cite{Berestycki15} studies
the necessary and sufficient condition for extinction.

Compared with the one-dimensional model \eqref{1.0}, the study of two-dimensional interaction of generalized nonlinear CSBP characterized by
the following SDE system is more challenging and the existing literature on this is sparse:
	\begin{equation}\label{1.1a}
		\left\{
		\begin{aligned}
			X_t &=X_0 +\int_0^t\theta_1(X_s,Y_s)\dd s+\int_0^t\gamma_{10}(X_s)\dd s+ \int_0^t\sqrt{\gamma_{11}(X_s)}\dd B_1(s) \\
			&~~~~
			+\int_0^t\int_0^\infty \int_0^{\gamma_{12}(X_{s-})}z\tilde{N}_1(\dd s,\dd z,\dd u),\\
			Y_t &=Y_0+\int_0^t \theta_2(Y_s,X_s)\dd s
			+\int_0^t\gamma_{20}(Y_s)\dd s +\int_0^t\sqrt{\gamma_{21}(Y_s)}\dd B_2(s) \\
			&~~~~
			+\int_0^t\int_0^\infty \int_0^{\gamma_{22}(Y_{s-})}z\tilde{N}_2(\dd s,\dd z,\dd u),
		\end{aligned}
		\right.
	\end{equation}
where $(B_1(t))_{t\ge0}$ and $(B_2(t))_{t\ge0}$ are Brownian motions,
and $\tilde{N}_1(\dd s,\dd z,\dd u)$ and $\tilde{N}_2(\dd s,\dd z,\dd u)$
are compensated Poisson random measures.
The SDE system \eqref{1.1a} is a solution to
the CSBP of two-types when $\gamma_{ij}(x)/x$ are nonnegative constants for all $i,j=1,2$.
If $\gamma_{i0}(x)=b_{i0}x$,
$\gamma_{i1}(x)=b_{i1}x$ with $b_{i1}>0$ and $\gamma_{i2}\equiv0$ for $i=1,2$,
then the SDE system \eqref{1.1a} is studied in \cite{Watanabe} and \cite{Cattiaux} when
$\theta_1(x,y)=a_1y,\theta_2(y,x)=a_2x$
and $\theta_1(x,y)=a_1xy,\theta_2(y,x)=a_2xy$
for $a_1,a_2>0$, respectively.
If $\gamma_{i0}(x)=a_{i0} x+\eta_i$,
$\gamma_{ij}(x)=b_{ij}x$ and $\theta_1(x,y)=a_1y,\theta_2(y,x)=a_2x$
with $b_{ij},a_i>0$ for $i,j=1,2$, then the SDEs system \eqref{1.1a} is studied
in \cite{Ma2013} and the more general system with multiple SDEs in \cite{bar}.
Two-type CSBP in varying environments characterized by the form of
\eqref{1.1a} is studied in \cite{LiZhanga,LiZhangb}.
The SDE system \eqref{1.1a} is also a
dynamic stochastic Lotka-Volterra-type population.
If $\gamma_{i0}(x)=b_{i0}x$ and $\gamma_{i1}(x)=b_{i1}x^2$ for $i=1,2$, $\theta_1(x,y)=a_2xy$, $\theta_2(y,x)=a_2xy$ with $b_{11},b_{21},a_1,a_2>0$,
$\gamma_{12}(x)=\gamma_{22}(x)\equiv0$, then
SDE \eqref{1.1a} is also called a competitive Lotka-Volterra model in random environments in \cite{Bao} and references, and \cite{Evans}
considers the case that $((B_1(t),B_2(t)))_{t\ge0}$ is a correlated two-dimensional
Brownian motion.
If $\gamma_{i2}(x)\equiv0$ and
$\gamma_{i1}(x)=b_{i1}x^2$ with $b_{i1}>0$ for $i=1,2$,
and $\theta_1(x,y)\theta_2(x,y)<0$,
then SDE \eqref{1.1a} is called stochastic predator-prey models, and we refer to \cite{DuDangYin} and \cite{BaoShao}  for related persistence/extinction results.
If $\theta_1(x,y)\equiv0$, $\theta_2(x,y)=:\theta(x)\kappa(y)<0$ and $\gamma_{10}(x),\gamma_{20}(y)\le0$, then the extinction-extinguishing dichotomy
for SDE system \eqref{1.1a} is discussed in \cite{RXYZ19}.
If $\theta_i(x,y)=\theta_i(x)\kappa_i(y)<0$ for $i=1,2$ and
$\gamma_{ij},\theta_i,\kappa_i$ are power functions for $i=1,2$ and $j=0,1,2$,
the SDE \eqref{1.1a} is mutually competitive and
rather sharp conditions on the extinction-extinguishing dichotomy are given in \cite{XYZh24}.
Our aim of this paper is to establish the extinction/non-extinction conditions for SDE \eqref{1.1a} with mutual enhancing,
that is, interacting functions $\theta_1,\theta_2>0$.

For simplicity and readability, in this paper we only consider the special form of \eqref{1.1a} with power function coefficients and power function intensities for $\tilde{N}_i(\dd s,\dd z,\dd u)$,
that is, the following SDE system:
	\begin{equation}\label{1.1}
		\left\{
		\begin{aligned}
			X_t &=X_0 +a_1\int_0^tX_s^{\theta_1}Y_s^{\kappa_1} \dd s-b_{10}\int_0^tX_s^{r_{10}}\dd s+\int_0^t\sqrt{2{\blue b_{11}}X_s^{r_{11}}}\dd B_1(s) \\
			&~~~~
			+\int_0^t\int_0^\infty \int_0^{b_{12}X_{s-}^{r_{12}}}z\tilde{N}_1(\dd s,\dd z,\dd u),\\
			Y_t &=Y_0+a_2\int_0^t Y_s^{\theta_2}X_s^{\kappa_2}\dd s
			-b_{20}\int_0^tY_s^{r_{20}}\dd s +\int_0^t\sqrt{2{\blue b_{21}}Y_s^{r_{21}}}\dd B_2(s) \\
			&~~~~
			+\int_0^t\int_0^\infty \int_0^{b_{22}Y_{s-}^{r_{22}}}z\tilde{N}_2(\dd s,\dd z,\dd u),
		\end{aligned}
		\right.
	\end{equation}
	where for $i=1,2,\,j=0,1,2$,  $a_i,\kappa_i>0$ and $\theta_i,r_{ij},b_{ij}\ge0$.
	For $i=1,2$, $(B_i(t))_{t\ge0}$ are two Brownian motions and
	$\tilde{N}_i(\dd s,\dd z,\dd u)$ are two compensated Poisson
	random measures with intensity $\dd s\mu_i(\dd z)\dd u$.
	Here
	$\mu_i(\dd z)=\frac{\alpha_i(\alpha_i-1)}{\Gamma(\alpha_i)\Gamma(2-\alpha_i)}z^{-1-\alpha_i}
	1_{\{z>0\}}\dd z$
	for $\alpha_i\in(1,2)$, $i=1,2$, and $\Gamma$ denotes the Gamma function.
	We always assume that $b_{11}+b_{12}>0$ and $b_{21}+b_{22}>0$.
	We also assume that
	$(B_1(t))_{t\ge0}$, $(B_2(t))_{t\ge0}$, $\{\tilde{N}_1(\dd s,\dd z,\dd u)\}$, and $\{\tilde{N}_2(\dd s,\dd z,\dd u)\}$ are independent of each other.

The main methods to develop the criteria for extinction/non-extinction
in \cite{MYZ21,baiyang} are adaptations of the approach for Chen's criteria on the uniqueness problem of Markov jump processes. These Chen's criteria are first established in \cite{Chen1986a,Chen1986b} and can also be found in
\cite[Theorems 2.25 and 2.27]{Chen04}. A similar approach to studying the boundary behaviors for Markov processes can also be found in \cite{MT93}.
Such an approach typically involves identifying an appropriate test function that is applied to the infinitesimal generator of the generalized nonlinear CSBP and proving the desired result  using a martingale argument.
In \cite{RXYZ19} a stochastic Lotka-Volterra-type population dynamical system is proposed as a solution to an SDE system with one-sided interaction.
The two-population model (\ref{1.1})  with mutually competitive interaction is further studied in \cite{XYZh24} where the coefficients $a_1$ and $a_2$ are both negative.
The main methods of \cite{RXYZ19} and \cite{XYZh24} rely on the generalization of Chen's  criteria technique from one-dimensional processes to two-dimensional processes.

In the paper, we consider the model of mutually enhancing populations, namely $a_1,a_2>0$, and leave the case of mixing interaction of $a_1a_2<0$ as a challenging {\em open} problem.
The main methods of this paper are to develop the criteria of extinction/non-extinction again
for a two-dimensional process which is an adaption of the approach for Chen's criteria again. The key to applying this criterion is to find appropriate test functions
for which our approach is mostly ad hoc, guessing bit by bit without an obvious  intuition. For the extinction behavior,
the criteria and the test functions in
\cite{RXYZ19,XYZh24} among other references are no longer applicable
in this paper, and we need to
develop new criteria and new testing functions that are
more demanding.

	For any generic stochastic  process $V:=(V(t))_{t\ge0}$ and constant $w>0$, let
	\beqnn
	\tau_0(V):=\inf\{t\ge0:V(t)=0\},\quad
	\tau_w^-(V):=\inf\{t\ge0:V(t)\le w\}
	\eeqnn
	and
	\beqnn
	\tau_w^+(V):=\inf\{t\ge0:V(t)\ge w\}
	\eeqnn
	with the convention $\inf\emptyset=\infty$.
	Let $\tau_0:=\tau_0(X)\wedge\tau_0(Y)$,
	$\tau_w^-:=\tau_w^-(X)\wedge\tau_w^-(Y)$
	and
	$\tau_w^+:=\tau_w^+(X)\wedge\tau_w^+(Y)$
	for $w>0$.
	In the following, we state the definition of a solution
	to the SDE system \eqref{1.1}, which is defined before the minimum of
	the first hitting time of 0 or the explosion time for the two processes
	$X$ and $Y$.
	
	\bdefinition\label{def}
	By a solution to SDE \eqref{1.1} we mean a two-dimensional c\`adl\`ag
	process $(X,Y):=((X_t,Y_t))_{t\ge0}$ satisfying
	SDE \eqref{1.1} up to $\gamma_n:=\tau_{1/n}^-
	\wedge\tau_n^+$ for each $n\ge1$ and
	$X_t=\limsup_{n\to\infty}X_{\gamma_n-}$
	and $Y_t=\limsup_{n\to\infty}Y_{\gamma_n-}$
	for $t\ge\lim_{n\to\infty}\gamma_n$.
	\edefinition

By Definition \ref{def},
$0$ and $\infty$ are absorbing states, and the solution is nonnegative.
Definition \ref{def} of the solution
allows for weaker conditions for the uniqueness of the solution. In particular, the existence and pathwise uniqueness of SDE \eqref{1.1} can be obtained by the same arguments as in \cite[Lemma A.1]{RXYZ19}.
Throughout this paper, we always assume that the  c\`adl\`ag process $(X,Y)$ is the unique solution to \eqref{1.1}, and consequently, the process $(X,Y)$ has the strong Markov property.
We
also assume that $X_0,Y_0>0$ are deterministic and that all stochastic processes are defined on the same filtered probability space $(\Omega,\mathscr{F},\mathscr{F}_t,\mathbf{P})$.
Let $\mathbf{E}$ denote the corresponding expectation.
For $i=1,2$, let
	\beqnn
	r_i
	 :=
	\min\big\{r_{ij}-\varrho_{ij},j\in\{b_{ij}\neq0,j=0,1,2\}\big\}, \quad
	b_i
	 :=
	\sum_{j=0,1,2}b_{ij}1_{\{r_i=r_{ij}-\varrho_{ij}\}},
	\eeqnn
where $\varrho_{ij}=j+1$ and $\varrho_{i2}=\alpha_i$ for $i=1,2$ and $j=0,1$.
	The above parameters determine the extinction behaviors for the
	processes $(X,Y)$.
	
	We first consider the non-extinction behaviors.
	Note that the nonlinear one-dimensional CSBP (with $a_1X_s^{\theta_1}Y_s^{\kappa_1}$
	replaced by $a_1X_s^{\theta_1}$ in \eqref{1.1})
 cannot reach $0$ if $r_1\ge(\theta_1-1)\wedge0$ (see \cite[p.2535]{LYZh}).
	The following result can also be shown similarly to  \cite[p.2535]{LYZh}.

	\bproposition\label{tp1.2}
	If $r_1,r_2\ge0$, then $\mbf{P}\{\tau_0<\infty\}=0$.
	\eproposition
	
	The next theorem concerns the extinction conditions
for the two-dimensional model \eqref{1.1} under which $Y_s^{\kappa_1}$
	in the interaction term
	$a_1X_s^{\theta_1}Y_s^{\kappa_1}$ does not affect the extinction of $X$
	when $r_2\ge0$ and $\theta_1-1<r_1<0$.
	
	\btheorem\label{t1.1a}
	Suppose that one of the following holds:
	\begin{itemize}
		\item[{\normalfont(i)}] $r_1\ge0$ and $\theta_2-1<r_2<0$;
		\item[{\normalfont(ii)}] $r_2\ge0$ and $\theta_1-1<r_1<0$.
	\end{itemize}
	Then $\mbf{P}\{\tau_0<\infty\}=0$.
	\etheorem

	The following theorem presents a non-extinction condition, when
		both $\theta_1-1<r_1<0$ and
		$\theta_2-1<r_2<0$ hold, which
	depends on the interaction term $a_1X_s^{\theta_1}Y_s^{\kappa_1}$
	and $a_2Y_s^{\theta_2}X_s^{\kappa_2}$.
	\btheorem\label{t1.1}
	$\mbf{P}\{\tau_0<\infty\}=0$ if $\theta_1-1<r_1<0$,
	$\theta_2-1<r_2<0$, and
	\begin{equation}\label{1.2}
		(r_1+1-\theta_1)(r_2+1-\theta_2)>\kappa_1\kappa_2.
	\end{equation}
	\etheorem

	Before considering  critical cases,
	we first state the following condition, which means
	that the drift term plays a dominant role
	when the processes are near zero.
	\bcondition\label{c1}
	\begin{itemize}
		\item[{\normalfont(i)}]
		For $b_{10}b_{20}\neq0$,
		{\blue\beqlb\label{1.8}
		r_{10}-1<\min\big\{r_{11}-\varrho_{1j}:
b_{1j}\neq0,j=1,2\big\}
		\eeqlb
		and
		\beqlb\label{1.9}
		r_{20}-1<\min\big\{r_{21}-\varrho_{2j}:
b_{2j}\neq0,j=1,2\big\},
		\eeqlb}
where $\varrho_{i1}=2$ and $\varrho_{i2}=\alpha_i$ for $i=1,2$;		
		\item[{\normalfont(ii)}]
		In addition to $(r_1+1-\theta_1)/\kappa_1<r_1/r_2$,
		either
		\eqref{1.8} holds and \eqref{1.9}
		is not satisfied for $b_{10}b_{20}\neq0$, or
		$b_{20}=0$ and
		\eqref{1.8} holds for $b_{10}\neq0$;
		\item[{\normalfont(iii)}]
		In addition to $(r_2+1-\theta_2)/\kappa_2<r_2/r_1$,
		either
		\eqref{1.9} holds and \eqref{1.8}
		is not satisfied for $b_{10}b_{20}\neq0$, or
		$b_{10}=0$ and
		\eqref{1.9} holds for $b_{20}\neq0$.
	\end{itemize}
	\econdition

	\btheorem\label{t1.4}
	Suppose that
	$\theta_1-1<r_1<0$,
	$\theta_2-1<r_2<0$, and
	\beqlb\label{1.2b}
	(r_1+1-\theta_1)(r_2+1-\theta_2)=\kappa_1\kappa_2.
	\eeqlb
	We also assume that one of the following holds:
	\begin{itemize}
		\item[{\normalfont(i)}]
		Condition \ref{c1}(i) is satisfied and
		\beqlb\label{1.2bb}
		(a_1/b_1)^{1/(r_1+1-\theta_1)}(a_2/b_2)^{1/\kappa_2}>1;
		\eeqlb
		\item[{\normalfont(ii)}] $b_{12}=b_{22}=0$,
		\beqnn
		r_{10}-1=r_{11}-2
		=r_{20}-1=r_{21}-2=:r,
		\eeqnn
		$\theta_1=\theta_2$, $\kappa_1=\kappa_2$,
		and $a_1a_2\ge b_1b_2$.
	\end{itemize}
	Then $\mbf{P}\{\tau_0<\infty\}=0$.
	\etheorem

Theorem \ref{t1.4}(ii) implies that  $\mbf{P}\{\tau_0<\infty\}=0$
		when the inequality in \eqref{1.2bb} is replaced by an equality.
Under \eqref{1.2b}, we have $(r_1+1-\theta_1)^{-1}=\frac{\kappa_1\kappa_2}{r_2+1-\theta_2}$
and then \eqref{1.2bb} is equivalent to
$(a_2/b_2)^{\frac{1}{r_2+1-\theta_2}}(a_1/b_1)^{\frac{1}{\kappa_1}}>1$,
which is a symmetric form of \eqref{1.2bb}. Next, we study the extinction behavior for $(X,Y)$.
	For the one-dimensional nonlinear CSBP (with $a_1X_s^{\theta_1}Y_s^{\kappa_1}$
	replaced by $a_1X_s^{\theta_1}$ in \eqref{1.1}),
	it reaches $0$ with a positive probability when $r_1<(\theta_1-1)\wedge0$.
	The following theorem shows that a similar result  still holds
	for the two-dimensional system \eqref{1.1} in this case.

	\btheorem\label{t1.2a}
	Suppose that
	one of the following holds:
	\begin{itemize}
		\item[{\normalfont(i)}]
		$r_1\le\theta_1-1$ and $r_1<0$;
		\item[{\normalfont(ii)}]
		$r_2\le\theta_2-1$ and $r_2<0$.
	\end{itemize}
	Then $\mbf{P}\{\tau_0<\infty\}>0$.
	\etheorem
	
		\bremark
		Observe that conditions (i) and (ii) of Theorem \ref{t1.2a} only concern processes $X$ and $Y$, respectively. This suggests that $\mbf{P}\{\tau_0(X)<\infty\}>0$ under condition (i) and $\mbf{P}\{\tau_0(Y)<\infty\}>0$ under condition (ii).
		\eremark


	In the following, we consider the case
	$\theta_1-1<r_1<0$ and $\theta_2-1<r_2<0$.
	We first state the following condition.

	\bcondition\label{c2}
	One of the following holds:
	\begin{itemize}
		\item[{\normalfont(i)}]
		Either
		\beqlb\label{3.113}
		\frac{r_1+1}{r_2+1}\vee\frac{r_1}{r_2}<\frac{\kappa_2}{r_2+1-\theta_2}
		\eeqlb
		or
		\beqlb\label{3.113b}
		\frac{r_2+1}{r_1+1}\vee\frac{r_2}{r_1}<\frac{\kappa_1}{r_1+1-\theta_1};
		\eeqlb
		\item[{\normalfont(ii)}]
Either
		\beqlb\label{3.122a}
		\frac{1-\theta_1}{\kappa_1-r_2}
		<\frac{\kappa_2}{r_2+1-\theta_2}\wedge
		\frac{\kappa_2-r_1}{1-\theta_2}\wedge(1-\theta_1)
		\wedge(2-\kappa_2)
		\eeqlb
or
		\beqlb\label{3.122b}
		\frac{1-\theta_2}{\kappa_2-r_1}
		<\frac{\kappa_1}{r_1+1-\theta_1}\wedge
		\frac{\kappa_1-r_2}{1-\theta_1}\wedge(1-\theta_2)
		\wedge(2-\kappa_1).
		\eeqlb
	\end{itemize}
	\econdition

The assumptions in Condition \ref{c2} are related to the choices of test functions (see \eqref{7.2} and Lemmas \ref{t7.4} and \ref{t7.5}) and may not allow intuitive interpretations.	
	
	\btheorem\label{t1.2}
	Suppose that $\theta_1-1<r_1<0$,
	$\theta_2-1<r_2<0$ and
	one of the following holds:
	\begin{itemize}
		\item[{\normalfont(i)}]
		Condition \ref{c1} holds and
		\beqlb\label{1.3}
		(r_1+1-\theta_1)(r_2+1-\theta_2)<\kappa_1\kappa_2;
		\eeqlb
		\item[{\normalfont(ii)}]
		For $i=1,2$,
		{\blue\beqlb\label{1.6}
				r_{i0}-1\ge\min\big\{r_{i1}-\varrho_{ij}:
b_{ij}\neq0,j=1,2\big\}, ~~~b_{i0}\neq0
		\eeqlb}
with $\varrho_{i1}=2$ and $\varrho_{i2}=\alpha_i$,		
		and \eqref{1.3} and Condition \ref{c2} holds;
		\item[{\normalfont(iii)}]
		For $i=1,2$, \eqref{1.6} holds,
		$r_1+1-\theta_1<\kappa_2$ and $r_2+1-\theta_2<\kappa_1$.
	\end{itemize}
	Then $\mbf{P}\{\tau_0<\infty\}>0$.
	\etheorem
	

	Theorem \ref{t1.2}(i) concerns the case in which drift terms dominate
when the process is near zero.
Theorems \ref{t1.2}(ii) and (iii) concern the case where diffusive terms dominate when the process is close to zero.
	When $r_1=r_2$, the inequality \eqref{1.3} implies Condition \ref{c2}(i).
In the following, we present a result
		when the  inequality in \eqref{1.3} is replaced by the equality.

	\btheorem\label{t1.4bb}
	Suppose that $\theta_1-1<r_1<0$, $\theta_2-1<r_2<0$
	and \eqref{1.2b} holds.
	We also assume that
	Condition \ref{c1} holds and
	\beqlb\label{1.3cc}
	(a_1/b_1)^{1/(r_1+1-\theta_1)}(a_2/b_2)^{1/\kappa_2}<1.
	\eeqlb
	Then $\mbf{P}\{\tau_0<\infty\}>0$.
	\etheorem

Given that \eqref{1.2b} holds, the inequality \eqref{1.3cc} is equivalent to
$(a_2/b_2)^{\frac{1}{r_2+1-\theta_2}}(a_1/b_1)^{\frac{1}{\kappa_1}}<1$,
which is a symmetric form of \eqref{1.3cc}.
It is more challenging and interesting to establish the
conditions for $\mbf{P}\{\tau_0<\infty\}=1$,
which will be considered in future work.
	Combining Proposition \ref{tp1.2}, Theorems \ref{t1.1a}--\ref{t1.1},
	\ref{t1.4}--\ref{t1.2a} and \ref{t1.2}--\ref{t1.4bb} we have the following remark.
	\bremark
	Suppose that
	$r_1=r_2:=r$.
	\begin{itemize}
		\item[{\normalfont(1)}]
		$\mbf{P}\{\tau_0<\infty\}=0$ if one of the following conditions:
		\begin{itemize}
			\item[{\normalfont(i)}] $r\ge0$;
			\item[{\normalfont(ii)}]
			$(\theta_1-1)\vee(\theta_2-1)<r<0$  and
			$(r+1-\theta_1)(r+1-\theta_2)>\kappa_1\kappa_2$.
		\end{itemize}
		\item[{\normalfont(2)}]
		$\mbf{P}\{\tau_0<\infty\}>0$ if one of the following conditions:
		\begin{itemize}
			\item[{\normalfont(i)}]
			$r\le[(\theta_1-1)\vee(\theta_2-1)]$ and $r<0$;
			\item[{\normalfont(ii)}]
			$(\theta_1-1)\vee(\theta_2-1)<r<0$, $(r+1-\theta_1)(r+1-\theta_2)<\kappa_1\kappa_2$, and
			in addition to either Condition  \ref{c1}(i) or \eqref{1.6} for $i=1,2$.
		\end{itemize}
		\item[{\normalfont(3)}]
		If Condition \ref{c1}(i) holds,
		$(\theta_1-1)\vee(\theta_2-1)<r<0$ and
		$(r+1-\theta_1)(r+1-\theta_2)=\kappa_1\kappa_2$,
		then
		\begin{itemize}
			\item[{\normalfont(i)}]
			$\mbf{P}\{\tau_0<\infty\}=0$ when
			$(a_1/b_1)^{1/(r+1-\theta_1)}(a_2/b_2)^{1/\kappa_2}>1$;
			\item[{\normalfont(ii)}]
			$\mbf{P}\{\tau_0<\infty\}>0$
			when
			$(a_1/b_1)^{1/(r+1-\theta_1)}(a_2/b_2)^{1/\kappa_2}<1$.
		\end{itemize}
	\end{itemize}
	\eremark

	At the end of this section, we give a brief description
of the approaches to proofs of theorems.
	We use the generalized Chen's criteria
	to establish the assertions on non-extinction and extinction with positive probability (see Section \ref{Criteria}).
	These Chen's criteria are first proposed in \cite{Chen1986a,Chen1986b} and can also be found in
	\cite[Theorems 2.25 and 2.27]{Chen04}.
	For the proofs of Proposition \ref{tp1.2},
Theorems \ref{t1.1a}--\ref{t1.1} and \ref{t1.4}, we apply the
	non-extinction criteria given in \cite{RXYZ19} (see, e.g., Proposition \ref{t3.1})
	to test functions that are  of
	two-variable power or logarithmic type
	(see the functions $g$ and $\tilde{g}$ defined in Lemmas \ref{t3.3} and \ref{t3.4}).
	For the proofs of Theorems \ref{t1.2a} and \ref{t1.2}--\ref{t1.4bb},
	we first establish the extinction criteria, which is given in Proposition \ref{t2.2}
	and a generalized Chen's criterion.
	The key test function for the proof of Theorem \ref{t1.2a}
	is chosen as a two-variable power function (see \eqref{6.2}).
	We select a two-variable power type test function in the proof of Theorems \ref{t1.2}--\ref{t1.4bb} (see \eqref{7.2}
	with the functions $h$ defined in \eqref{3.7} and \eqref{3.7b} and Lemmas \ref{t7.4}, \ref{t7.5}
and \ref{t7.7}
	for different conditions).
	For the case where the drift term  dominates (Condition \ref{c1})
	as the system approaches zero, the test functions $h$ and $\tilde{h}$ are given by \eqref{3.7} and \eqref{3.7b}, respectively,
	and condition (iii) in Proposition \ref{t2.2} is verified
	by Lemmas \ref{t7.2} and \ref{t7.3}.
	For the case where diffusive terms dominate (see \eqref{1.6})
and Condition \ref{c2},
	we modify the function $h$ defined in \eqref{3.7}.
	Under Condition \ref{c2}(i), we first replace $y^{\rho_2}$ with $y-y^{1-\varepsilon}$ for a sufficiently small $y,\varepsilon>0$ in \eqref{3.7} and then change it into $\tilde{h}(y)$ given by \eqref{7.5}
	(see \eqref{7.3}).
	For Condition \ref{c2}(ii),
	we also change $\rho_1\ge1$ to $0<\rho_1<1$ and smooth it
	by replacing $x^{\rho_1}$ with $(x+y^\delta)^{\rho_1}$ for some $\delta>1$
	(see \eqref{7.4}).
	For condition (iii) in Theorem \ref{t1.2},
	the function $h$ is given in Lemma \ref{t7.7}.
	The key condition (iii) in Proposition \ref{t2.2} is verified by
	Lemmas \ref{t7.4}, \ref{t7.5}, and \ref{t7.7}.

	Similarly, the results can be proved for the cases of $b_{10},b_{20}<0$
 and $b_{10}b_{20}<0$, and we omit them
	here.
	Let $C^2((0,\infty))$ and $C^2((0,\infty)\times(0,\infty))$
	denote the second-order continuous differentiable function spaces on $(0,\infty)$
	and $(0,\infty)\times(0,\infty)$, respectively.
	The remainder of the paper
	is arranged as follows.
	We state the extinction criteria for the two-dimensional
	process in Section \ref{Criteria}.
	The proofs for the assertions in this section
	are given in Section 3.

	\section{Criteria for extinction}\label{Criteria}
	\setcounter{equation}{0}
	
	In this section, we establish some criteria
that will be used to prove the theorems in Section \ref{intro}.
	In the following, let $((x_t,y_t))_{t\ge0}$ be a two-dimensional c\`adl\`ag
	process with deterministic $x_0,y_0>0$, where $(x_t)_{t\ge0}$ and $(y_t)_{t\ge0}$ are two nonnegative
	processes defined before the minimum of their first times of hitting
	zero or explosion.
	Let $\mathcal{L}$ denote the operator such that for each $g\in C^2((0,\infty)\times(0,\infty))$,
	the process
	\beqlb\label{3.1}
	t\mapsto M^g_{t\wedge\gamma_{m,n}} \quad \mbox{is a martingale},
	\eeqlb
	where
	\beqnn
	M^g_t:=g(x_t,y_t)-g(x_0,y_0)-\int_0^t\mathcal{L}g(x_s,y_s)\dd s
	\eeqnn
	and $\gamma_{m,n}:=\tau_{1/m}^-\wedge\tau_n^+$
	with $\tau_{1/m}^-:=\tau^-_{1/m}(x)\wedge\tau^-_{1/m}(y)$
	and $\tau^+_n:=\tau^+_n(x)\wedge\tau^+_n(y)$.
	In this section, let $\tau_0:=\tau_0(x)\wedge\tau_0(y)$
	and $\tau_\infty:=\lim_{n\to\infty}\tau_n^{+}$.
	The following two criteria on non-extinction and extinction of process
	$((x_t,y_t))_{t\ge0}$, which generalize Chen's criteria for the uniqueness problem of Markov jump processes. For the SDE system, the operator $\mathcal{L}$ can be obtained by It\^o's formula.

The following Proposition \ref{t3.1} given in \cite{RXYZ19} is a criterion
 to establish the non-extinction assertion and test functions $g$ often take either the form of  two-dimensional power type functions with negative power or the form of logarithm type functions, which is a key point.
This proposition is used to prove Proposition \ref{tp1.2},
Theorems \ref{t1.1a}--\ref{t1.1} and \ref{t1.4}.

	\bproposition\label{t3.1}
	(\cite[Proposition 2.1]{RXYZ19})
	For any fixed $n\ge1$,
	if there is a nonnegative function $g\in C^2((0,\infty)\times(0,\infty))$
	and a constant $d_n>0$ such that the following hold:
	(i) $\lim_{x\wedge y\to0}g(x,y)=\infty$;
	(ii) $\mathcal{L}g(x,y)\le d_ng(x,y)$ for all $0<x,y\le n$,
	then $\mbf{P}\{\tau_0<\infty\}=0$.
	\eproposition

The proposition \ref{t2.2} gives a new criterion for
extinction and can be shown by adapting the approach to Chen's criteria.
The key to applying Proposition \ref{t2.2} is to select a suitable test function
and to verify the key condition (iii) in the following.
This proposition is used to prove Theorems \ref{t1.2a}, \ref{t1.2} and  \ref{t1.4bb}.
The test functions $g$ are often of the form of two-dimensional power type functions with positive powers, which are mostly obtained by ad hoc  guessing  without much intuitive interpretation, and the test functions in the previous papers are no longer appropriate.
Note that the criteria and  test functions  in Theorems \ref{t1.2a}, \ref{t1.2}, and  \ref{t1.4bb} are different from those
in \cite{RXYZ19,XYZh24} and other references.
The proof of Proposition \ref{t2.2} is based on
\eqref{3.1} and a martingale argument.

	\bproposition\label{t2.2}
	Let $v>0$ and $0<x_0,y_0<v$ be fixed.
	Suppose that there are a function $g\in C^2((0,\infty)\times(0,\infty))$
	and a constant $d>0$ such that
	\begin{itemize}
		\item[{\normalfont(i)}]
		$g(x_0,y_0)>0$ and $C_0:=\sup_{x,y>0}g(x,y)<\infty$;
		\item[{\normalfont(ii)}]
		$g(x,y)\le0$ for all $(x,y)$ with $x\vee y\geq v$;
		\item[{\normalfont(iii)}]
		$\mathcal{L}g(x,y)\ge dg(x,y)$ for all $0<x,y\le v$.
	\end{itemize}
	Then $\mbf{P}\{\tau_0<\infty\}\ge g(x_0,y_0)/C_0>0$.
	\eproposition
	\proof
	In view of \eqref{3.1}, for all large enough $m\ge1$,
	\beqnn
	\mbf{E}\big[g(x_{t\wedge\gamma_{m,v}},y_{t\wedge\gamma_{m,v}})\big]
	=
	g(x_0,y_0)+\int_0^t\mbf{E}\big[
	\mathcal{L}g(x_s,y_s)1_{\{s\le\gamma_{m,v}\}}\big]\dd s.
	\eeqnn
	It then follows from integration by parts that
	\beqnn
	\ar\ar
	\e^{-dt}\mbf{E}\big[g(x_{t\wedge\gamma_{m,v}},y_{t\wedge\gamma_{m,v}})\big] \cr
	\ar=\ar
	g(x_0,y_0)+\int_0^t\e^{-ds}
	\dd\Big(\mbf{E}\big[g(x_{s\wedge\gamma_{m,v}},y_{s\wedge\gamma_{m,v}})\big]\Big)
	+\int_0^t
	\mbf{E}\big[g(x_{s\wedge\gamma_{m,n}},y_{s\wedge\gamma_{m,v}})\big]
	\dd\e^{-ds} \cr
	\ar=\ar
	g(x_0,y_0)
	+\int_0^t\e^{-ds} \mbf{E}\big[
	\mathcal{L}g(x_s,y_s)1_{\{s\le\gamma_{m,v}\}}\big] \dd s
	-d\int_0^t
	\e^{-ds} \mbf{E}\big[g(x_{s\wedge\gamma_{m,v}},y_{s\wedge\gamma_{m,v}})\big]
	\dd s.
	\eeqnn
	Under conditions (i) and (iii), we have
	\beqnn
	\ar\ar
	C_0\e^{-dt}-g(x_0,y_0)\ge
	\e^{-dt}\mbf{E}\big[g(x_{t\wedge\gamma_{m,v}},y_{t\wedge\gamma_{m,v}})\big]
	-g(x_0,y_0) \cr
	\ar\ge\ar
	d\int_0^t\e^{-ds} \mbf{E}\big[
	g(x_s,y_s)1_{\{s\le\gamma_{m,v}\}}\big] \dd s
	-d\int_0^t
	\e^{-ds} \mbf{E}\big[g(x_{s\wedge\gamma_{m,v}},y_{s\wedge\gamma_{m,v}})\big]
	\dd s \cr
	\ar=\ar
	-d\int_0^t\e^{-ds} \mbf{E}\big[
	g(x_{s\wedge\gamma_{m,v}},y_{s\wedge\gamma_{m,v}})1_{\{s>\gamma_{m,v}\}}\big] \dd s \cr
	\ar\ge\ar
	-d\int_0^t\e^{-ds} \mbf{E}\big[
	g^+(x_{s\wedge\gamma_{m,v}},y_{s\wedge\gamma_{m,v}})1_{\{s>\gamma_{m,v}\}}\big] \dd s,
	\eeqnn
	where $g^+(u):=g(u)\vee 0$.
	Letting $t\to\infty$ we obtain
	\beqnn
	g(x_0,y_0)
	\ar\le\ar
	d\int_0^\infty\e^{-ds} \mbf{E}\big[
	g^+(x_{s\wedge\gamma_{m,v}},y_{s\wedge\gamma_{m,v}})1_{\{s>\gamma_{m,v}\}}\big] \dd s \cr
	\ar=\ar
	\mbf{E}\Big[g^+(x_{\gamma_{m,v}},y_{\gamma_{m,v}})
	d\int_{\gamma_{m,v}}^\infty\e^{-ds}\dd s\Big]
	=\mbf{E}\big[g^+(x_{\gamma_{m,v}},y_{\gamma_{m,v}})
	\e^{-d\gamma_{m,v}}\big].
	\eeqnn
	Then, by the dominated convergence, we get
	\beqlb\label{2.2}
	g(x_0,y_0)
	\ar\le\ar
	\lim_{m\to\infty}\mbf{E}\big[g^+(x_{\gamma_{m,v}},y_{\gamma_{m,v}})
	\e^{-d\gamma_{m,v}}\big]
	=\mbf{E}\big[g^+(x_{\tau_0\wedge\tau^+_v},y_{\tau_0\wedge\tau^+_v})
	\e^{-d(\tau_0\wedge\tau^+_v)}\big] \cr
	\ar=\ar
	\mbf{E}\Big[g^+(x_{\tau_0\wedge\tau^+_v},y_{\tau_0\wedge\tau^+_v})
	\e^{-d(\tau_0\wedge\tau^+_v)}
	\big(1_{\{\tau_0<\tau^+_v\}}+1_{\{\tau_0\ge\tau^+_v,\tau^+_v<\infty\}}
	+1_{\{\tau_0=\tau^+_v=\infty\}}\big)\Big] \cr
	\ar=\ar
	\mbf{E}\Big[g^+(x_{\tau_0},y_{\tau_0})
	\e^{-d\tau_0}
	1_{\{\tau_0<\tau^+_v\}}
	+
	g^+(x_{\tau^+_v},y_{\tau^+_v})
	\e^{-d\tau^+_v}1_{\{\tau_0\ge\tau^+_v,\tau^+_v<\infty\}}\Big]  \cr
	\ar\le\ar
	C_0\mbf{P}\{\tau_0<\tau^+_v\},
	\eeqlb
	where conditions (i) and (ii) are used in the last inequality.
	Since $g(x_0,y_0)>0$ under condition (i),
	then by \eqref{2.2},
	\beqnn
	\mbf{P}\{\tau_0<\infty\}\ge\mbf{P}\{\tau_0<\tau^+_v\}\ge
	g(x_0,y_0)/C_0>0,
	\eeqnn
	which completes the proof.
	\qed

	\section{Proofs of the main results}
	\setcounter{equation}{0}
	
In preparation for the main proofs,
	we first introduce some notation and inequalities.
	For any $g\in C^2((0,\infty)\times(0,\infty))$ and $x,y,z\ge0$ define
	\beqlb\label{3.21}
	K_z^1g(x,y):=
	g(x+z,y)-g(x,y)-g'_x(x,y)z
	\eeqlb
	and
	\beqlb\label{3.22}
	K_z^2g(x,y):=
	g(x,y+z)-g(x,y)-g'_y(x,y)z.
	\eeqlb
By Taylor's formula, for any bounded function $g$ with continuous second derivative,
 \beqlb\label{6.0a}
g(x+z)-g(x)-g'(x)z=z^2\int_0^1g''(x+zu)(1-u)\dd u.
 \eeqlb
For the Borel function $g$ in $(0,\infty)$, replacing the variable $z$ with $zx$ we get
 \beqlb\label{6.0b}
\int_0^\infty g(z)\mu_i(\dd z)=
x^{-\alpha_i}\int_0^\infty g(zx)\mu_i(\dd z),\qquad i=1,2,
 \eeqlb
where we use the fact $\mu_i(\dd (zx))=x^{-\alpha_i}\mu_i(\dd z)$.
The operator $\mathcal{L}$ is given by
	\beqlb\label{3.0}
	\mathcal{L}g(x,y)
	\ar:=\ar
	a_1x^{\theta_1}y^{\kappa_1}g'_x(x,y)
	+a_2y^{\theta_2}x^{\kappa_2}g'_y(x,y) \cr
	\ar\ar
	-b_{10}x^{r_{10}}g'_x(x,y)
	+ b_{11}x^{r_{11}}g''_{xx}(x,y)
	+b_{12}x^{r_{12}}\int_0^\infty K_z^1g(x,y)\mu_1(\dd z) \cr
	\ar\ar
	-b_{20}y^{r_{20}}g'_y(x,y)
	+ b_{21}y^{r_{21}}g''_{yy}(x,y)
	+b_{22}y^{r_{22}}\int_0^\infty K_z^2g(x,y)\mu_2(\dd z).
	\eeqlb
	For $i=1,2$ and $\rho<\alpha_i$ let $c_i(\rho):=\frac{\Gamma(\alpha_i-\rho)}{
		\Gamma(\alpha_i)\Gamma(2-\rho)}$.
	By \cite[Lemma 4.2]{XYZh24} and  \eqref{6.0a},
	for $\rho\in(-\infty,0)\cup(0,1)\cup(1,\alpha_i)$,
	\beqlb\label{3.0a}
	c_i(\rho)=[\rho(\rho-1)]^{-1}\int_0^\infty [(1+z)^\rho-1-\rho z]\mu_i(\dd z)
	=\int_0^\infty z^2\mu_i(\dd z)\int_0^1(1+zu)^{\rho-2}(1-u)\dd u.
	\eeqlb
	We now present a lemma for the proofs below.
	\blemma\label{t2.3}
	\begin{itemize}
		\item[{\normalfont(i)}]
		For any $u,v\ge0$ and $p,q>1$ with $1/p+1/q=1$, we have
		$u+v\ge p^{1/p}q^{1/q} u^{1/p}v^{1/q}$ and
		$u/p+v/q\ge u^{1/p}v^{1/q}$.
		\item[{\normalfont(ii)}]
		For $x,y\ge0$, we have
		$x^p+y^p\ge (x+y)^p$ for any $0<p\le 1$ and
		$x^p+y^p\ge 2^{1-p}(x+y)^p$  for any $p>1$.
		\item[{\normalfont(iii)}]
		Suppose that $p_1,p_2,p_3,p_4>0$ and $c_1,c_2,c_3>0$.
		If $p_3/p_1+p_4/p_2>1$,
		then there is a constant $0<c<1$ such that $c_1x^{p_1}+c_2y^{p_2}
		\ge c_3x^{p_3}y^{p_4}$ for all $0<x,y<c$.
	\end{itemize}
	\elemma
	\proof
	Assertion (i) follows immediately from the Young inequality, and
	assertion (ii) is obvious.
	Let $0<\delta<1$ satisfy
	$1-p_4/p_2<\delta<p_3/p_1$.
	Then $\delta p_1-p_3<0$ and $(1-\delta) p_2-p_4<0$.
	By assertion (i), there are constants $c_4,c>0$ such that
	\beqnn
	c_1x^{p_1}+c_2y^{p_2}\ge
	c_4x^{\delta p_1}y^{(1-\delta)p_1}
	=c_3x^{p_3}y^{p_4}
	(c_4/c_3)x^{\delta p_1-p_3}y^{(1-\delta)p_1-p_4}
	\ge
	c_3x^{p_3}y^{p_4}
	\eeqnn
	for all $0<x,y<c$.
	This proves assertion (iii).
	\qed

	\subsection{Proofs of Proposition \ref{tp1.2} and Theorems \ref{t1.1a}--\ref{t1.1} and \ref{t1.4}}

In this subsection, we use Proposition \ref{t3.1}
to complete the proofs of Proposition \ref{tp1.2} and  Theorems \ref{t1.1a}--\ref{t1.1} and \ref{t1.4}.
The test functions $g$ are presented in Lemmas \ref{t3.3} and \ref{t3.4}.
The key to the proofs is to verify condition (ii) in Proposition \ref{t3.1} for
which we need  some lemmas in the following subsubsection.

	\subsubsection{Preliminaries}

	The following lemma is used to prove Lemma \ref{t3.3} in the next subsection where the key to its proof is applying Lemma \ref{t2.3}(i).
	
	\blemma\label{t3.2}
	Suppose that $r_2>\theta_2-1$, $\rho_1,\rho_2>\theta_1\vee\theta_2$ and
	\beqlb\label{5.20}
	\frac{r_2+1-\theta_2}{\kappa_2}
	>
	\frac{1+\rho_2+\kappa_1-\theta_2}{1+\rho_1+\kappa_2-\theta_1}.
	\eeqlb
	Then there are constants $\delta,C>0$ such that
	\beqnn
	y^{\theta_2-1-\rho_2}x^{\kappa_2}
	+x^{\theta_1-1-\rho_1}y^{\kappa_1}
	\ge Cm^{\delta}y^{r_2-\rho_2}
	\eeqnn
	\elemma
	for all $m\ge1$ and $0<x,y\le m^{-1}$.
	
	\proof
	Let
	\beqnn
	p:=
	\frac{1+\rho_1+\kappa_2-\theta_1}{1+\rho_1-\theta_1},~~
	q:=\frac{1+\rho_1+\kappa_2-\theta_1}{\kappa_2}.
	\eeqnn
	Then
	$p,q>1$ and $p^{-1}+q^{-1}=1$.
	Moreover, by \eqref{5.20}, we obtain
	\beqnn
	\frac{1}{q}
	=
	\frac{\kappa_2}{1+\rho_1+\kappa_2-\theta_1}
	<
	\frac{r_2+1-\theta_2}{1+\rho_2+\kappa_1-\theta_2},
	\eeqnn
	and then
	\beqnn
	\delta
	:=
	(r_2+1-\theta_2)-\frac{1+\rho_2+\kappa_1-\theta_2}{q}>0,~
	\frac{1+\rho_1+\kappa_2-\theta_1}{p}-(1+\rho_1-\theta_1)=0.
	\eeqnn
	Now, by Lemma \ref{t2.3}(i),
	\beqnn
	x^{1+\rho_1+\kappa_2-\theta_1}+y^{1+\rho_2+\kappa_1-\theta_2}
	\ge
	p^{1/p}q^{1/q}
	x^{(1+\rho_1+{\blue\kappa_2}-\theta_1)/p}
	y^{(1+\rho_2+\kappa_1-\theta_2)/q},
	\eeqnn
	and then
	\beqnn
	(x^{1+\rho_1+\kappa_2-\theta_1}+y^{1+\rho_2+\kappa_1-\theta_2})
	/(y^{r_2+1-\theta_2}x^{1+\rho_1-\theta_1})
	\ge
	p^{1/p}q^{1/q}
	y^{-\delta}.
	\eeqnn
	It follows that for $m\ge1$ and $0<x,y\le m^{-1}$,
	\beqnn
	\ar\ar
	y^{\theta_2-1-\rho_2}x^{\kappa_2}+x^{\theta_1-1-\rho_1}y^{\kappa_1} \cr
	\ar=\ar
	y^{r_2-\rho_2}(x^{1+\rho_1+\kappa_2-\theta_1}+y^{1+\rho_2+\kappa_1-\theta_2})
	/(y^{r_2+1-\theta_2}x^{1+\rho_1-\theta_1}) \cr
	\ar\ge\ar
	p^{1/p}q^{1/q}
	y^{r_2-\rho_2}y^{-\delta}
	\ge
	p^{1/p}q^{1/q}m^{\delta}
	y^{r_2-\rho_2},
	\eeqnn
	which completes the proof. \qed

The following Lemmas \ref{t5.3} and \ref{t5.1} are used to establish
Lemma \ref{t3.4} in the next subsubsection and the key to the proofs is to apply
Lemma \ref{t2.3}(i) again.
	For $\theta_1-1<r_1$ and
	$\theta_2-1<r_2$,
	under \eqref{1.2b}, there are constants $\rho_1,\rho_2>\theta_1\vee\theta_2$  such that
	\beqlb\label{1.22b}
	\frac{r_1+1-\theta_1}{\kappa_1}
	= \frac{\kappa_2}{r_2+1-\theta_2}
	=
	\frac{1+\rho_1+\kappa_2-\theta_1}{1+\rho_2+\kappa_1-\theta_2}.
	\eeqlb
	
	\blemma\label{t5.3} Recalling $a_i$ and $b_i$ defined in Section \ref{intro},
	suppose that $r_1>\theta_1-1$,
	$r_2>\theta_2-1$ and
	\beqlb\label{5.2}
	(a_1/b_1)^{1/(r_1+1-\theta_1)}(a_2/b_2)^{1/\kappa_2}\ge1,
	\eeqlb
	and that  \eqref{1.22b} holds for  $\rho_1,\rho_2>\theta_1\vee\theta_2$.
	Then there is a constant $\delta_0>0$ such that
	\beqnn
	a_1\rho_1\delta_0 x^{\theta_1-1-\rho_1}y^{\kappa_1}
	+a_2\rho_2x^{\kappa_2}y^{\theta_2-1-\rho_2}
	\ge
	b_1\rho_1 \delta_0 x^{r_1-\rho_1}
	+
	b_2\rho_2y^{r_2-\rho_2},\quad x,y>0.
	\eeqnn
	\elemma
	\proof
	Let
	\beqnn
	\ar\ar
	p_1:=
	\frac{1+\rho_2+\kappa_1-\theta_2}{1+\rho_2-\theta_2},\quad
	q_1:=\frac{1+\rho_2+\kappa_1-\theta_2}{\kappa_1}, \cr
	\ar\ar
	p_2:=\frac{1+\rho_1+\kappa_2-\theta_1}{\kappa_2},\quad
	q_2:=\frac{1+\rho_1+\kappa_2-\theta_1}{1+\rho_1-\theta_1}.
	\eeqnn
	Then $p_1,q_1,p_2,q_2>1$, $p_1^{-1}+q_1^{-1}=1$ and
	$p_2^{-1}+q_2^{-1}=1$.
	In view of \eqref{1.22b}, we obtain
	$q_1=(1+\rho_1+\kappa_2-\theta_1)/(r_1+1-\theta_1)$.
	Combining this with \eqref{5.2} we have
	$(b_1/a_1)^{q_1}(b_2/a_2)^{p_2}\le1$,
	which implies that
	\beqnn
	(b_1\rho_1)^{q_1}(b_2\rho_2)^{p_2}
	\le
	(a_1\rho_1)^{q_1}(a_2\rho_2)^{p_2}
	=
	(a_1\rho_1)^{q_1/p_1}(a_2\rho_2)^{p_2/q_2}
	a_1\rho_1 a_2\rho_2
	\eeqnn
	and then
	\beqnn
	\frac{(b_2\rho_2)^{p_2}}{(a_2\rho_2)^{p_2/q_2}a_1\rho_1 p_2/q_1}
	\le
	\frac{(a_1\rho_1 )^{q_1/p_1}a_2\rho_2q_1/p_2}{(b_1\rho_1)^{q_1}}.
	\eeqnn
	Thus there is a constant $\delta_0>0$ satisfying
	\beqnn
	\frac{(b_2\rho_2)^{p_2}}{(a_2\rho_2)^{p_2/q_2}a_1\rho_1 p_2/q_1}
	\le\delta_0\le
	\frac{(a_1\rho_1 )^{q_1/p_1}a_2\rho_2q_1/p_2}{(b_1\rho_1)^{q_1}}.
	\eeqnn
	It follows that
	\beqlb\label{5.13}
	(a_1\rho_1\delta_0)^{1/p_1}
	(a_2\rho_2q_1/p_2)^{1/q_1}\ge b_1\rho_1 \delta_0,\quad
	(a_1\rho_1\delta_0 p_2/q_1)^{1/p_2}
	(a_2\rho_2)^{1/q_2}
	\ge b_2\rho_2.
	\eeqlb

	By \eqref{1.22b} again,
	\beqnn
	\frac{1}{q_1}
	\ar=\ar
	\frac{\kappa_1}{1+\rho_2+\kappa_1-\theta_2}
	=
	\frac{r_1+1-\theta_1}{1+\rho_1+\kappa_2-\theta_1}, \cr
	\frac{1}{p_2}
	\ar=\ar
	\frac{\kappa_2}{1+\rho_1+\kappa_2-\theta_1}
	=
	\frac{r_2+1-\theta_2}{1+\rho_2+\kappa_1-\theta_2},
	\eeqnn
	and then
	\beqnn
	\ar\ar
	(r_1+1-\theta_1)-\frac{1+\rho_1+\kappa_2-\theta_1}{q_1}=0,\quad
	\frac{1+\rho_2+\kappa_1-\theta_2}{p_1}-(1+\rho_2-\theta_2)=0, \cr
	\ar\ar
	(r_2+1-\theta_2)-\frac{1+\rho_2+\kappa_1-\theta_2}{p_2}=0,\quad
	\frac{1+\rho_1+\kappa_2-\theta_1}{q_2}-(1+\rho_1-\theta_1)=0.
	\eeqnn
	Now by Lemma \ref{t2.3}(i) and \eqref{5.13},
	\beqnn
	\ar\ar
	p_1^{-1}a_1\rho_1\delta_0 y^{1+\rho_2+\kappa_1-\theta_2}
	+p_2^{-1}a_2\rho_2x^{1+\rho_1+\kappa_2-\theta_1} \cr
	\ar\ge\ar
	p_1^{1/p_1}q_1^{1/q_1}
	(p_1^{-1}a_1\rho_1\delta_0 y^{1+\rho_2+\kappa_1-\theta_2})^{1/p_1}
	(p_2^{-1}a_2\rho_2x^{1+\rho_1+\kappa_2-\theta_1})^{1/q_1} \cr
	\ar=\ar
	(a_1\rho_1 \delta_0)^{1/p_1}(a_2\rho_2q_1/p_2)^{1/q_1}
	x^{{\blue r_1+1}-\theta_1}y^{1+\rho_2-\theta_2}
	\ge
	b_1\rho_1 \delta_0 x^{{\blue r_1+1}-\theta_1}y^{1+\rho_2-\theta_2}
	\eeqnn
	and
	\beqnn
	\ar\ar
	q_1^{-1}a_1\rho_1\delta_0 y^{1+\rho_2+\kappa_1-\theta_2}
	+q_2^{-1}a_2\rho_2x^{1+\rho_1+\kappa_2-\theta_1} \cr
	\ar\ge\ar
	p_2^{1/p_2}q_2^{1/q_2}
	(q_1^{-1}a_1\rho_1\delta_0 y^{1+\rho_2+\kappa_1-\theta_2})^{1/p_2}
	(q_2^{-1}a_2\rho_2x^{1+\rho_1+\kappa_2-\theta_1})^{1/q_2} \cr
	\ar=\ar
	(a_1\rho_1\delta_0 p_2/q_1)^{1/p_2}
	(a_2\rho_2)^{1/q_2}y^{{\blue r_2+1}-\theta_2}x^{1+\rho_1-\theta_1}
	\ge
	b_2\rho_2y^{{\blue r_2+1}-\theta_2}x^{1+\rho_1-\theta_1}.
	\eeqnn
	Then one can conclude the proof.
	\qed
	
	Similar to Lemma \ref{t5.3} we have the following lemma.
	\blemma\label{t5.1}
	If both the assumption of Theorem \ref{t1.4}(ii) and
	\eqref{1.2b} hold,
	and $r>\theta-1$, then
	there is a constant $\delta_0>0$ such that
	\beqnn
	\delta_0 a_1 y^{1+\kappa-\theta}+a_2 x^{1+\kappa-\theta}
	\ge
	\delta_0 b_1x^{\kappa}y^{1-\theta}
	+
	b_2y^{\kappa}x^{1-\theta},\qquad x,y>0.
	\eeqnn
	\elemma
	\proof
	Let $p:=(1+\kappa-\theta)/(1-\theta)$ and $q:=p /(p-1)$.
	As the essential same argument in the proof of Lemma \ref{t5.3},
	we obtain
	\beqnn
	p^{-1}\delta_0 a_1x^{\theta-1}y^{\kappa}
	+p^{-1}a_2x^{\kappa}y^{\theta-1}
	\ge
	\delta_0 b_1x^{\kappa+\theta-1}, ~~
	q^{-1}\delta_0 a_1y^{\kappa}x^{\theta-1}
	+q^{-1}a_2y^{\theta-1}x^{\kappa}
	\ge
	b_2 y^{\kappa+\theta-1}
	\eeqnn
	for some $\delta_0>0$,
	which implies the assertion.
	\qed

	\subsubsection{Proofs of Proposition \ref{tp1.2} and  Theorems \ref{t1.1a}--\ref{t1.1} and \ref{t1.4}}
	
	We next apply Proposition \ref{t3.1} to complete the proofs of
	Theorems \ref{t1.1a}--\ref{t1.1} and \ref{t1.4}, where the following
Lemmas \ref{t3.3} and \ref{t3.4} are applied to validate  condition (ii) of Proposition \ref{t3.1}.
	
	\blemma\label{t3.3}
	Given $g(x,y):=x^{-\rho_1}+y^{-\rho_2}, \, x,y>0$,
	if the assumptions of either Theorem \ref{t1.1a} or \ref{t1.1} hold,
	then there are constants $\rho_1,\rho_2,C_n>0$ such that
	$\mathcal{L}g(x,y)\le C_ng(x,y)$ for all $0<x,y\le n$ and $n\ge1$.
	\elemma
	\proof
	(i) We first prove the result under the assumptions of Theorem \ref{t1.1a}.
	Since the proofs are similar, we only present that under condition (i)
	of Theorem \ref{t1.1a}.
	Observe that there are constants $\rho_1,\rho_2>\theta_1\vee\theta_2$
	such that \eqref{5.20} holds.
	
	By \eqref{6.0b} and \eqref{3.0a},
	\beqlb\label{3.3}
	\int_0^\infty [(x+z)^{-\rho_i}-x^{-\rho_i}+\rho_i zx^{-\rho_i}]\mu_i(\dd z)
	=c_i\rho_i(\rho_i+1)x^{-\alpha_i-\rho_i}
	\eeqlb
	with $c_i:=c_i(-\rho_i)$.
	Then by \eqref{3.0}, for all $0<x,y\le n$,
	\beqlb\label{5.6}
	\mathcal{L}g(x,y)
	\ar=\ar
	-a_1\rho_1x^{\theta_1-1-\rho_1}y^{\kappa_1}
	-a_2\rho_2y^{\theta_2-1-\rho_2}x^{\kappa_2}  \cr
	\ar\ar
	+\rho_1 b_{10}x^{r_{10}-1-\rho_1}
	+b_{11} \rho_1(\rho_1+1)x^{r_{11}-2-\rho_1}
	+ c_1\rho_1(\rho_1+1)b_{12}x^{r_{12}-\alpha_1-\rho_1} \cr
	\ar\ar
	+\rho_2 b_{20}y^{r_{20}-1-\rho_2}
	+b_{21} \rho_2(\rho_2+1)y^{r_{21}-2-\rho_2}
	+ c_2\rho_2(\rho_2+1)b_{22}y^{r_{22}-\alpha_2-\rho_2} \cr
	\ar\le\ar
	-a_1\rho_1x^{\theta_1-1-\rho_1}y^{\kappa_1}
	-a_2\rho_2y^{\theta_2-1-\rho_2}x^{\kappa_2} \cr
	\ar\ar
	+\rho_1(\rho_1+1)x^{-\rho_1}\big[b_{10}x^{r_{10}-1}
	+b_{11} x^{r_{11}-2}+b_{12}c_1 x^{r_{12}-\alpha_1} \big] \cr
	\ar\ar
	+\rho_2(\rho_2+1)y^{-\rho_2}\big[b_{20}y^{r_{20}-1}
	+b_{21} y^{r_{21}-2}
	+b_{22}c_2 y^{r_{22}-\alpha_2} \big] \cr
	\ar\le\ar
	-a_1\rho_1x^{\theta_1-1-\rho_1}y^{\kappa_1}
	-a_2\rho_2y^{\theta_2-1-\rho_2}x^{\kappa_2}
	+C_{1,n}x^{r_1-\rho_1}
	+C_{2,n}y^{r_2-\rho_2}
	\eeqlb
	with $
	C_{i,n}:=\rho_i(\rho_i+1)
	[b_{i0}n^{r_{i0}-1-r_i}
	+b_{i1} n^{r_{i1}-2-r_i}
	+b_{i2} n^{r_{i2}-\alpha_i-r_i}]$
	for $i=1,2$.
	By \eqref{5.20} and Lemma \ref{t3.2},
	\beqlb\label{5.7}
	\ar\ar
	-a_1\rho_1x^{\theta_1-1-\rho_1}y^{\kappa_1}
	-a_2\rho_2y^{\theta_2-1-\rho_2}x^{\kappa_2}
	+C_{2,n}y^{r_2-\rho_2} \cr
	\ar\le\ar
	-[(a_1\rho_1)\wedge(a_2\rho_2)]
	\cdot[x^{\theta_1-1-\rho_1}y^{\kappa_1}+y^{\theta_2-1-\rho_2}x^{\kappa_2}]
	+C_{2,n}y^{r_2-\rho_2} \cr
	\ar\le\ar
	\big[C_{2,n}-C[(a_1\rho_1)\wedge(a_2\rho_2)]
	m^{\delta}\big]y^{r_2-\rho_2}\le0,\quad 0<x,y\le m^{-1}
	\eeqlb
	for all large enough $m$, where $\delta,C>0$ are the constants
	determined in Lemma \ref{t3.2}.
	From \eqref{5.6} and \eqref{5.7} it follows that
	\beqlb\label{5.8}
	\mathcal{L}g(x,y)
	\ar\le\ar
	-a_1\rho_1x^{\theta_1-1-\rho_1}y^{\kappa_1}
	-a_2\rho_2y^{\theta_2-1-\rho_2}x^{\kappa_2}
	+C_{2,n}y^{r_2-\rho_2}
	+C_{1,n}n^{r_1}g(x,y) \cr
	\ar\le\ar
	C_{1,n}n^{r_1}g(x,y),\qquad 0<x,y\le m^{-1}
	\eeqlb
	for all large enough $m$.
	Under the assumptions in condition (i), by \eqref{5.6}, we have
	\beqlb\label{5.9}
	\mathcal{L}g(x,y)
	\le
	C_{1,n}x^{r_1-\rho_1}
	+C_{2,n}y^{r_2-\rho_2}
	\le
	[C_{1,n}n^{r_1}
	+C_{2,n}m^{-r_2}]g(x,y)
	\eeqlb
	for all $m^{-1}\le y\le n$ and $0<x\le n$,
	and
	\beqlb\label{5.11}
	\mathcal{L}g(x,y)
	\ar\le\ar
	-a_2\rho_2y^{\theta_2-1-\rho_2}x^{\kappa_2}
	+C_{1,n}x^{r_1-\rho_1}
	+C_{2,n}y^{r_2-\rho_2} \cr
	\ar\le\ar
	y^{r_2-\rho_2}[
	C_{2,n}-a_2\rho_2y^{\theta_2-1-r_2}x^{\kappa_2}]
	+C_{1,n}n^{r_1}g(x,y) \cr
	\ar\le\ar
	y^{r_2-\rho_2}[
	C_{2,n}-a_2\rho_2m_0^{r_2+1-\theta_2}m^{-\kappa_2}]
	+C_{1,n}n^{r_1}g(x,y) \cr
	\ar\le\ar
	C_{1,n}n^{r_1}g(x,y),\qquad
	m^{-1}\le x\le n,~0<y\le m_0^{-1}
	\eeqlb
	for $m\ge1$ and large enough $m_0\ge1$. Observe that
	$\mathcal{L}g(x,y)\le C_ng(x,y)$ for some $C_n>0$ and all
	$m^{-1}\le x\le n$ and $m_0^{-1}\le y\le n$.
	Combining \eqref{5.8}--\eqref{5.11} we get the assertion.
	
	(ii) Now we show the result under the assumptions
	of Theorem \ref{t1.1}.
	Under \eqref{1.2},
	there are constants
	$\rho_1,\rho_2>\theta_1\vee\theta_2$ such that
	\beqnn
	\frac{{\blue r_1+1}-\theta_1}{\kappa_1}
	>
	\frac{1+\rho_1+\kappa_2-\theta_1}{1+\rho_2+\kappa_1-\theta_2}
	> \frac{\kappa_2}{{\blue r_2+1}-\theta_2}.
	\eeqnn
	Then by Lemma \ref{t3.2}
	and  the same arguments as in \eqref{5.7},
	we have
	\beqnn
	-a_1\rho_1x^{\theta_1-1-\rho_1}y^{\kappa_1}
	-a_2\rho_2y^{\theta_2-1-\rho_2}x^{\kappa_2}
	+C_{1,n}x^{r_1-\rho_1}
	+C_{2,n}y^{r_2-\rho_2}\le0,~~0<x,y\le m^{-1}
	\eeqnn
	for large enough $m\ge1$.
	By  arguments similar to those in \eqref{5.9} and \eqref{5.11},
	we have $\mathcal{L}g(x,y)\le C_{n,m}g(x,y)$ for some constant $C_{n,m}>0$
	for large enough $m\ge1$ when $m^{-1}\le x\le n$ and $0<y\le n$,
	or $m^{-1}\le y\le n$ and $0<x\le n$.
	Now $\mathcal{L}g(x,y)\le C_{n,m}g(x,y)$ for $0<x,y\le n$
	by \eqref{5.6}.
	\qed

	\blemma\label{t3.4}
	Suppose that the assumptions of Theorem \ref{t1.4} hold.
	Under the conditions of Theorem \ref{t1.4}(i) let
	\beqnn
	\tilde{g}(x,y):=\delta_0x^{-\rho_1}+y^{-\rho_2},\qquad x,y>0,
	\eeqnn
	and under the conditions of Theorem \ref{t1.4}(ii)
	define a nonnegative function $\tilde{g}\in C^2((0,\infty)\times(0,\infty))$ such that
	\beqnn
	\tilde{g}(x,y):=(\delta_0+1) \ln n
	+\delta_0 \ln x^{-1}+\ln y^{-1},\qquad 0<x,y<n.
	\eeqnn
	Then there are constants $\rho_1,\rho_2,\delta_0,C_n>0$ such that
	$\mathcal{L}\tilde{g}(x,y)\le C_n\tilde{g}(x,y)$ for all $0<x,y<n$ and $n\ge1$.
	\elemma
	\proof
	(i) We first prove it under the conditions in Theorem \ref{t1.4}(i).
	Under \eqref{1.2b},
	there are constants
	$\rho_1,\rho_2>\theta_1\vee\theta_2$ such that
	\eqref{1.22b} holds.
	Let $\delta_0>0$ be the constant determined in the following.
	Then by \eqref{3.0} and \eqref{3.3},
	we obtain
	$\mathcal{L}\tilde{g}(x,y)
	=
	-J_1(x,y)+J_2(x,y)+J_3(x,y)$
	with $J_1(x,y)
	:=
	\delta_0 a_1\rho_1x^{\theta_1-1-\rho_1}y^{\kappa_1}
	+a_2\rho_2y^{\theta_2-1-\rho_2}x^{\kappa_2}$,
	\beqnn
	J_2(x,y)
	:=
	\delta_0 b_{10} \rho_1x^{r_{10}-1-\rho_1}
	+b_{20}\rho_2y^{r_{20}-1-\rho_2}
	=
	\delta_0 b_1 \rho_1x^{r_1-\rho_1}
	+b_2\rho_2y^{r_2-\rho_2}
	\eeqnn
	under Condition \ref{c1}(i) and
	\beqnn
	J_3(x,y)
	\ar:=\ar
	\delta_0 b_{11} \rho_1(\rho_1+1)x^{r_{11}-2-\rho_1}
	+b_{12}\delta_0 c_1(-\rho_1)\rho_1(\rho_1+1) x^{r_{12}-\alpha_1-\rho_1}\cr
	\ar\ar
	+b_{21}\rho_2(\rho_2+1)y^{r_{21}-2-\rho_2}
	+b_{22}c_2(-\rho_2)\rho_2(\rho_2+1) y^{r_{22}-\alpha_2-\rho_2} \cr
	\ar\le\ar
	c_0[x^{r_{11}-2}+x^{r_{12}-\alpha_1}]x^{-\rho_1}
	+c_0[y^{r_{21}-2}+y^{r_{22}-\alpha_2}]y^{-\rho_2}
	\eeqnn
	for some constant $c_0>0$.
	Under \eqref{1.2bb},
	there is a small enough constant $\varepsilon>0$ such that
	\beqnn
	\Big[\frac{a_1}{(1+\varepsilon)b_1}\Big]^{1/(r_1+1-\theta_1)}
	\Big[\frac{a_2}{(1+\varepsilon)b_2}\Big]^{1/\kappa_2}\ge1.
	\eeqnn
	By Lemma \ref{t5.3}, there is a constant $\delta_0>0$ such that
	$(1+\varepsilon)J_2(x,y)\le J_1(x,y)$ for all $x,y>0$.
	By Condition \ref{c1}(i) again, there is a constant $c_n>0$
	such that
	\beqnn
	-\varepsilon\delta_0 b_{10} \rho_1x^{r_{10}-1}
	+c_0[x^{r_{11}-2}+x^{r_{12}-\alpha_1}]\le c_n\delta_0,
	~
	-\varepsilon b_{20}\rho_2y^{r_{20}-1}+c_0[y^{r_{21}-2}+y^{r_{22}-\alpha_2}]
	\le c_n
	\eeqnn
	for all $0<x,y\le n$.
	Thus, for all $0<x,y\le n$, $J_3(x,y)-\varepsilon J_2(x,y)
	\le c_n\tilde{g}(x,y)$ and we have
	$\mathcal{L}\tilde{g}(x,y)\le c_n\tilde{g}(x,y)$.
	
	(ii)
	We then prove the result under the conditions in Theorem \ref{t1.4}(ii).
	Let $\theta:=\theta_1$ and $\kappa:=\kappa_1$.
Then $r=\kappa+\theta-1$.
	By Lemma \ref{t5.1}, there is a constant $\delta_0>0$ such that
	\beqnn
	I(x,y):=\delta_0 a_1 y^{1+\kappa-\theta}+a_2 x^{1+\kappa-\theta}
	-\delta_0 b_1x^{\kappa}y^{1-\theta}-b_2y^{\kappa}x^{1-\theta}\ge0,
	\qquad x,y>0.
	\eeqnn
	It thus follows from \eqref{3.0} that
	\beqnn
	\mathcal{L}\tilde{g}(x,y)
	\ar=\ar
	-\delta_0 a_1x^{\theta-1} y^\kappa-a_2 y^{\theta-1} x^{\kappa}
	+\delta_0(b_{10}x^{r_{10}-1}+b_{11}x^{r_{11}-2}) \cr
	\ar\ar
	+(b_{20}y^{r_{20}-1}+b_{21}y^{r_{21}-2})
	=
	-x^{\theta-1}y^{\theta-1}I(x,y)\le0,
	\eeqnn
	which ends the proof.
	\qed
	
Now we are ready to complete the proofs of Proposition \ref{tp1.2} and Theorems \ref{t1.1a}, \ref{t1.1}
and \ref{t1.4}.

\noindent{\it Proof of Proposition \ref{tp1.2}.}
Let $g$ be the test function given  in Lemma \ref{t3.3}.
By \eqref{5.6},
 \beqnn
\mathcal{L}g(x,y)
\le C_{1,n}x^{r_1-\rho_1}
	+C_{2,n}y^{r_2-\rho_2}
\le [C_{1,n}x^{r_1}+C_{2,n}y^{r_2}]g(x,y)
\le
[C_{1,n}n^{r_1}+C_{2,n}n^{r_2}]g(x,y)
 \eeqnn
for some constants $C_{1,n},C_{2,n}\ge0$
and all $0<x,y\le n$. The conclusion follows from Proposition \ref{t3.1}. \qed
	
	\noindent{\it Proof of Theorems \ref{t1.1a} and \ref{t1.1}.}
	Let $g$ be the test function determined in Lemma \ref{t3.3}.
	Then by Proposition \ref{t3.1}  and Lemma \ref{t3.3} one completes the proofs.
	\qed

	\noindent{\it Proof of Theorem \ref{t1.4}.}
	Let $\tilde{g}$ be the test function defined in Lemma \ref{t3.4}.
	Then one completes the proof combining Proposition \ref{t3.1} and Lemma \ref{t3.4}.
	\qed

	\subsection{Proofs of Theorems \ref{t1.2a}, \ref{t1.2} and  \ref{t1.4bb}}

In this subsection, we use Proposition \ref{t2.2} to establish the proofs of Theorems \ref{t1.2a}, \ref{t1.2} and  \ref{t1.4bb}.
For the proof of Theorem \ref{t1.2a},
the test function $g$ is taken as  a power function with positive power
(see \eqref{6.2}).
For the proofs of Theorems \ref{t1.2} and \ref{t1.4bb},
the test functions $g(x,y):=c+h(x,y)$ for some constant $c>0$
(see \eqref{7.2}) and $h$ are modifications of $-(x^{\rho_1}+y^{\rho_2})^\rho$ in different cases.
The key to applying Proposition \ref{t2.2} is verifying assumption (iii) of Proposition \ref{t2.2}.
For the case of dominant drifts (Condition \ref{c1}),
the function $h$ is chosen as $h(x,y):=-(x^{\rho_1}+y^{\rho_2})^\rho$
for $\rho_1,\rho_2>1$ and $0<\rho<1$ and the estimation of $\mathcal{L}h$ is the key to verifying
assumption (iii) of Proposition \ref{t2.2},
which are given in Subsubsection \ref{Equations {1.1} with dominate drifts}.
For the case of dominant diffusion or jump terms
(assumption \eqref{1.6}), the functions $h$ are modifications of $-(x^{\rho_1}+y)^\rho$
or $-(x+y^{\rho_2})^\rho$ for $\rho_1,\rho_2>0$ and $0<\rho<1$
and the estimation of $\mathcal{L}h$ is also the key in checking
assumption (iii) of Proposition \ref{t2.2}, which are given in Subsubsection
\ref{Equations {1.1} with dominate diffusive or jump terms}.
In Subsubsection \ref{Proof of Theorem {t1.2a}}, the proof of Theorem \ref{t1.2a}
is presented. Three lemmas  in Subsubsection \ref{Preliminaries} are used to estimate $\mathcal{L}h$.
The proofs of Theorems \ref{t1.2} and  \ref{t1.4bb} are given
at the end of Subsubsection
\ref{Equations {1.1} with dominate diffusive or jump terms}.

\subsubsection{Proof of Theorem \ref{t1.2a}}\label{Proof of Theorem {t1.2a}}	

In this subsubsection, we prove  Theorem \ref{t1.2a} by
applying Proposition \ref{t2.2} to the test function $g$  given by \eqref{6.2}.
	
	\noindent{\it Proof of Theorem \ref{t1.2a}.}
	Since the proofs are similar, we only present that under condition (i).
	Let $0<\rho<1\wedge(-r_1)$ be fixed as follows.
	Since $r_1\le \theta_1-1$ and $\kappa_1>0$, there are constants $0<v<1$
	and $d>0$ such that
	\beqlb\label{6.1}
	2^{-1}b_1(1-\rho)[1\wedge c_1(\rho)]  x^{r_1+1-\theta_1}
	-a_1 y^{\kappa_1}\ge0,~
	2^{-1}b_1\rho(1-\rho)[1\wedge c_1(\rho)]
	x^{\rho+r_1}-a_2\ge dv
	\eeqlb
	for all $0<x,y\le v$, where $c_1(\rho)$ is given by \eqref{3.0a}.
	Let
	\beqlb\label{6.2}
	g(x,y):=v-x^{\rho}-y,\qquad x,y>0.
	\eeqlb
	Then $\sup_{x,y>0}g(x,y)\le v$, and $g(x,y)\le0$
	for all $x\ge v$ or $y\ge v$.
	By \eqref{3.0a},
	\beqnn
	x^{-\rho+\alpha_1}\int_0^\infty [(x+z)^\rho-x^\rho-\rho zx^{\rho-1}]\mu_1(\dd z)
	=\int_0^\infty [(1+z)^\rho-1-\rho z]\mu_1(\dd z)
	=\rho(\rho-1)c_1(\rho).
	\eeqnn
 Then it follows from \eqref{3.0}  and \eqref{6.1}  that
	\beqnn
	\mathcal{L}g(x,y)
	\ar=\ar
	-a_1\rho x^{\rho+\theta_1-1}y^{\kappa_1}-a_2y^{\theta_2}x^{\kappa_2}
	+b_{20} y^{r_{20}}
	+b_{10}\rho x^{\rho+r_{10}-1} \cr
	\ar\ar
	+b_{11}\rho(1-\rho)x^{\rho+r_{11}-2}
	+b_{12}c_1(\rho)\rho(1-\rho)x^{\rho+r_{12}-\alpha_1} \cr
	\ar\ge\ar
	-a_1\rho x^{\rho+\theta_1-1}y^{\kappa_1}-a_2y^{\theta_2}x^{\kappa_2}
	+b_1\rho(1-\rho)[1\wedge c_1(\rho)]x^{\rho+r_1} \cr
	\ar\ge\ar
	2^{-1}b_1\rho(1-\rho)[1\wedge c_1(\rho)]x^{\rho+r_1}-a_2 \cr
	\ar\ar
	+\rho x^{\rho+\theta_1-1}\Big[2^{-1}b_1(1-\rho)[1\wedge c_1(\rho)] x^{r_1+1-\theta_1}
	-a_1 y^{\kappa_1}\Big] \cr
	\ar\ge\ar
	dv\ge d g(x,y),\qquad  0<x,y\le v.
	\eeqnn
	Using Proposition \ref{t2.2} one gets $\mbf{P}\{\tau_0<\infty\}\ge
g(X_0,Y_0)/v$ for $0<X_0,Y_0<(2^{-1}v)^{1/\rho}$.
In general, $X_0,Y_0>0$,  we conclude the proof by applying the strong Markov property and the same arguments
	at the end of the proof of \cite[Theorem 1.3]{XYZh24}.
	\qed

	\subsubsection{Preliminaries}\label{Preliminaries}
	
	The following Lemmas \ref{t7.1}--\ref{t4.1a} are needed
to estimate $\mathcal{L}h$ in Lemmas \ref{t7.2} and \ref{t7.3}
with $h$ defined by \eqref{3.7} and \eqref{3.7b}. Lemma \ref{t2.3}(iii) and (i) are applied to show Lemmas \ref{t7.1} and \ref{t5.2}, respectively.  Lemma \ref{t4.1a} contains an estimate of the integral $\int_0^\infty K_z^i h(x,y)\mu_i(\dd z)$ for $i=1,2$ in \eqref{3.0}, which is  applied to prove Lemma \ref{t7.4}.

	\blemma\label{t7.1}
	Suppose that
	\beqlb\label{3.6}
	\frac{r_1+1-\theta_1}{\kappa_1}
	<
	\frac{r_1+\rho_1}{r_2+\rho_2}
	<\frac{\kappa_2}{r_2+1-\theta_2}
	\eeqlb
   for $\theta_1-1< r_1<0$, $\theta_2-1< r_2<0$ and  $\rho_1,\rho_2\ge1$.
	Then for any constants $c_1,c_2,c_3,c_4>0$ there is a constant $0<c<1$ such that for all $0<x,y\le c$,
	\beqnn
	\ar\ar
	c_1x^{r_1+\rho_1}+c_2y^{r_2+\rho_2}
	-c_3x^{\theta_1-1+\rho_1}y^{\kappa_1}
	-c_4y^{\theta_2-1+\rho_2}x^{\kappa_2}\ge0.
	\eeqnn
	\elemma
	\proof
	Under \eqref{3.6}, we have
	\beqnn
	\frac{\theta_1-1+\rho_1}{r_1+\rho_1}
	+\frac{\kappa_1}{r_2+\rho_2}>1,\quad
	\frac{\kappa_2}{r_1+\rho_1}
	+\frac{\theta_2-1+\rho_2}{r_2+\rho_2}>1.
	\eeqnn
	Thus, by Lemma \ref{t2.3}(iii), there is a constant $0<c<1$ such that
	\beqnn
	2^{-1}c_1x^{r_1+\rho_1}+2^{-1}c_2y^{r_2+\rho_2}
	\ge c_3x^{\theta_1-1+\rho_1}y^{\kappa_1},\quad
	2^{-1}c_1x^{r_1+\rho_1}+2^{-1}c_2y^{r_2+\rho_2}
	\ge c_4y^{\theta_2-1+\rho_2}x^{\kappa_2}
	\eeqnn
	for all $0<x,y\le c$.
	This concludes the proof.
	\qed

	\blemma\label{t5.2}
	Suppose that
	\beqlb\label{4.0a}
	\Big(\frac{a_1}{(1-\varepsilon_0)b_1}\Big)^{1/(r_1+1-\theta_1)}
	\Big(\frac{a_2}{(1-\varepsilon_0)b_2}\Big)^{1/\kappa_2}<1
	\eeqlb
	for {\blue $\theta_1-1< r_1$, $\theta_2-1< r_2$}
		and $\varepsilon_0\in(0,1)$,
	and in addition,
	\beqlb\label{4.0}
	\big(a_1/\tilde{b}_1\big)^{1/(r_1+1-\theta_1)}
	\big(a_2/\tilde{b}_2\big)^{1/\kappa_2}<1
	\eeqlb
	and
	\beqlb\label{5.1}
	\frac{r_1+1-\theta_1}{\kappa_1}
	=
	\frac{r_1+\rho_1}{r_2+\rho_2}
	=\frac{\kappa_2}{r_2+1-\theta_2},
	\eeqlb
	for $\rho_1,\rho_2>1$, $0<\rho<\rho_1^{-1}\wedge\rho_2^{-1}$ and $\tilde{b}_i:=(1-\varepsilon_0)(1-\rho_i\rho)b_i$, $i=1,2$.
	Then there is $\delta_0>0$ such that
	\beqnn
	\delta_0 \tilde{b}_1\rho_1x^{r_1+\rho_1}+
	\tilde{b}_2\rho_2y^{r_2+\rho_2}
	>
	\delta_0 a_1\rho_1 x^{\theta_1+\rho_1-1}y^{\kappa_1}
	+a_2\rho_2 y^{\theta_2+\rho_2-1}x^{\kappa_2}.
	\eeqnn
{\blue When \eqref{5.1} holds and $(a_1/b_1\big)^{1/(r_1+1-\theta_1)}
	(a_2/b_2\big)^{1/\kappa_2}=1$, 
there is $\delta_0>0$ such that
	\beqnn
	\delta_0 b_1\rho_1x^{r_1+\rho_1}+
	b_2\rho_2y^{r_2+\rho_2}
	\ge
	\delta_0 a_1\rho_1 x^{\theta_1+\rho_1-1}y^{\kappa_1}
	+a_2\rho_2 y^{\theta_2+\rho_2-1}x^{\kappa_2}.
	\eeqnn}
	\elemma
	\proof
{\blue Since the proofs are similar, we only state the first one.}
 Let
	\beqnn
	p_1:=\frac{r_1+\rho_1}{\theta_1-1+\rho_1},~ q_1:=\frac{r_2+\rho_2}{\kappa_1},\quad
	p_2:=\frac{r_1+\rho_1}{\kappa_2},~
	q_2:=\frac{r_2+\rho_2}{\theta_2-1+\rho_2}.
	\eeqnn
   By \eqref{5.1} we have
	\beqlb\label{5.3}
	\frac{1}{q_1}=\frac{\kappa_1}{r_2+\rho_2}=\frac{r_1+1-\theta_1}{r_1+\rho_1},
	\quad
	\frac{1}{p_2}=\frac{\kappa_2}{r_1+\rho_1}=\frac{r_2+1-\theta_2}{r_2+\rho_2}.
	\eeqlb
	It follows that
	\beqlb\label{5.21}
	p_1^{-1}+q_1^{-1}=1,\quad p_2^{-1}+q_2^{-1}=1,
	\eeqlb
which implies that $p_1,q_1,p_2,q_2>1$.

	Under \eqref{4.0} and \eqref{5.1},
	we have
	\beqnn
	1
	\ar>\ar
	[(a_1/\tilde{b}_1)^{1/(r_1+1-\theta_1)}
	(a_2/\tilde{b}_2)^{1/\kappa_2}]^{r_1+\rho_1} \cr
	\ar=\ar
	(a_1/\tilde{b}_1)^{(r_2+\rho_2)/\kappa_1}
	(a_2/\tilde{b}_2)^{(r_1+\rho_1)/\kappa_2}
	=
	a_1^{q_1}a_2^{p_2}/(\tilde{b}_1^{q_1}\tilde{b}_2^{p_2}),
	\eeqnn
	where \eqref{5.3} is used in the last equality.
	Then
	\beqnn
	(a_1\rho_1)^{q_1}(a_2\rho_2)^{p_2}
	<
	(\tilde{b}_1\rho_1)^{q_1}(\tilde{b}_2\rho_2)^{p_2}
	=(\tilde{b}_1\rho_1)^{q_1/p_1}
	\tilde{b}_2\rho_2
	\cdot(\tilde{b}_2\rho_2)^{p_2/q_2} \tilde{b}_1\rho_1,
	\eeqnn
	which implies that
	\beqnn
	\frac{(a_2\rho_2)^{p_2}}{
		(\tilde{b}_2\rho_2)^{p_2/q_2} \tilde{b}_1\rho_1p_2/q_1}<
	\frac{(\tilde{b}_1\rho_1)^{q_1/p_1}
		\tilde{b}_2\rho_2q_1/p_2}{(a_1\rho_1)^{q_1}}.
	\eeqnn
	Then there is a constant $\delta_0>0$ such that that
	\beqnn
	\frac{(a_2\rho_2)^{p_2}}{(
		\tilde{b}_2\rho_2)^{p_2/q_2} \tilde{b}_1\rho_1p_2/q_1}
	<\delta_0<
	\frac{(\tilde{b}_1\rho_1)^{q_1/p_1}
		\tilde{b}_2\rho_2q_1/p_2}{(a_1\rho_1)^{q_1}}.
	\eeqnn
	Therefore,
	\beqnn
	\ar\ar
	(q_1/p_2)^{1/q_1}(\delta_0 \tilde{b}_1\rho_1)^{1/p_1}
	(\tilde{b}_2\rho_2)^{1/q_1}>\delta_0 a_1\rho_1,\quad
	(p_2/q_1)^{1/p_2}(\delta_0 \tilde{b}_1\rho_1)^{1/p_2}
	(\tilde{b}_2\rho_2)^{1/q_2}>a_2\rho_2.
	\eeqnn
	Now by Lemma \ref{t2.3}(i) and \eqref{5.3} again
	\beqnn
	\ar\ar
	p_1^{-1}\delta_0 \tilde{b}_1\rho_1 x^{r_1+\rho_1}+
	p_2^{-1}\tilde{b}_2\rho_2y^{r_2+\rho_2} \cr
	\ar\ge\ar
	p_1^{1/p_1}q_1^{1/q_1}
	(p_1^{-1}\delta_0 \tilde{b}_1\rho_1 x^{r_1+\rho_1})^{1/p_1}
	(p_2^{-1}\tilde{b}_2\rho_2y^{r_2+\rho_2})^{1/q_1} \cr
	\ar=\ar
	(q_1/p_2)^{1/q_1}(\delta_0 \tilde{b}_1\rho_1)^{1/p_1}
	(\tilde{b}_2\rho_2)^{1/q_1}x^{\theta_1+\rho_1-1}y^{\kappa_1}
	>
	\delta_0 a_1\rho_1
	x^{\theta_1+\rho_1-1}y^{\kappa_1}
	\eeqnn
	and
	\beqnn
	\ar\ar
	q_1^{-1}\delta_0 \tilde{b}_1\rho_1x^{r_1+\rho_1}+
	q_2^{-1}\tilde{b}_2\rho_2y^{r_2+\rho_2} \cr
	\ar\ge\ar
	p_2^{1/p_2}q_2^{1/q_2}
	(q_1^{-1}\delta_0 \tilde{b}_1\rho_1 x^{r_1+\rho_1})^{1/p_2}
	(q_2^{-1}\tilde{b}_2\rho_2y^{r_2+\rho_2})^{1/q_2} \cr
	\ar=\ar
	(p_2/q_1)^{1/p_2}(\delta_0 \tilde{b}_1\rho_1)^{1/p_2}
	(b_2\rho_2)^{1/q_2}y^{\theta_2+\rho_2-1}x^{\kappa_2}
	>
	a_2\rho_2 y^{\theta_2+\rho_2-1}x^{\kappa_2}.
	\eeqnn
	This completes the proof.
	\qed

To estimate the integral $\int_0^\infty K_z^i h(x,y)\mu_i(\dd z)$
for $i=1,2$ in \eqref{3.0} with the function $h$ defined in
\eqref{3.7}, \eqref{3.7b} and in Lemmas \ref{t7.4} and \ref{t7.5}
for different cases, we introduce functions $K(v,z)$ and $M(x,y)$ in the following. Given $0<\rho<1$ and $\rho_1\ge1$, for $0\le v\le 1$ and $z>0$ let
	\beqlb\label{3.7bb}
	K(v,z):=-\big(v[(1+z)^{\rho_1}-1]+1\big)^\rho+1+zv\rho\rho_1.
	\eeqlb
Then for $M(x,y):=-(x^{\rho_1}+y)^\rho$,
we have
$K_{zx}^1M(x,y)=(x^{\rho_1}+y)^\rho
K(\frac{x^{\rho_1}}{x^{\rho_1}+y},z)$.
It thus follows from \eqref{6.0b} that
 \beqlb\label{3.10}
\int_0^\infty K_z^1M(x,y)\mu_1(\dd z)
 \ar=\ar
x^{-\alpha_1}
\int_0^\infty K_{zx}^1M(x,y)\mu_1(\dd z) \cr
 \ar=\ar
x^{-\alpha_1}(x^{\rho_1}+y)^\rho
\int_0^\infty K(\frac{x^{\rho_1}}{x^{\rho_1}+y},z)\mu_1(\dd z).
 \eeqlb
In the following we give an estimate on $\int_0^\infty K(v,z)\mu_1(\dd z)$.
The estimate mainly relies on Taylor's formula in \eqref{6.0a}
and change of variable similar to \eqref{6.0b}.
	
	\blemma\label{t4.1a}
	(i) There is a constant $\rho_0\ge2$ such that for all $\rho_1>\rho_0$,
all $\tilde{\rho}\in(0,1)$ and $\rho:=\tilde{\rho}/\rho_1$,
	\beqlb\label{3.28}
	\int_0^\infty K(v,z)\mu_1(\dd z)
	\ge
	\rho(1-\rho)\rho_1^2 v^2 d_1
	-\rho\rho_1(\rho_1-1) v  \tilde{d}_1
	\eeqlb
	for some constants $d_1,\tilde{d}_1>0$.

	(ii) For any $\rho_1>1$ and $0<\rho\rho_1<1$, and for all $\delta>0$ we have
	\beqlb\label{3.11a}
	\int_0^\infty K(v,z)\mu_1(\dd z)\ge
	\rho\rho_1(1-\rho\rho_1) c_1(\rho\rho_1)v^2
	-\rho\rho_1(\rho_1-1)[v(1-v)d_{1,\delta}+d_{2,\delta}v^\rho],
	\eeqlb
	where
	$d_{1,\delta}:=\int_0^\delta z^2(1+z)^{\rho_1-2}\mu_1(\dd z)$,
	$d_{2,\delta}:=\int_\delta^\infty z^2 \mu_1(\dd z)\int_0^1(1+uz)^{\rho\rho_1-2}\dd u$ and function $c_1$ is defined in (\ref{3.0a}).
	\elemma
	\proof
	(i)
	Observe that
	\beqlb\label{3.8a}
	K''_{zz}(v,z)
	\ar=\ar
	\rho(1-\rho)\rho_1^2v^2(1+z)^{2\rho_1-2}\big[v[(1+z)^{\rho_1}-1]+1\big]^{\rho-2} \cr
	\ar\ar
	-\rho\rho_1(\rho_1-1)v(1+z)^{\rho_1-2}\big[v[(1+z)^{\rho_1}-1]+1\big]^{\rho-1} \cr
	\ar=:\ar
	\rho(1-\rho)\rho_1^2 H_1(v,z)-\rho\rho_1(\rho_1-1)H_2(v,z).
	\eeqlb
	Then $H_1(v,z)\ge v^2(1+z)^{\rho\rho_1-2}$
	and $H_2(v,z)\le v(1+z)^{\rho_1-2}$ for $0\le v\le1$.
	Now by \eqref{6.0a},
	\beqlb\label{3.8}
	\int_0^1 K(v,z)\mu_1(\dd z)
	\ar=\ar
	\int_0^1 z^2\mu_1(\dd z)\int_0^1
	K''_{zz}(v,zu) (1-u)\dd u \cr
	\ar\ge\ar
	\rho(1-\rho)\rho_1^2v^2\int_0^1 z^2\mu_1(\dd z)\int_0^1
	(1+zu)^{\rho\rho_1-2}(1-u)\dd u \cr
	\ar\ar
	-\rho\rho_1(\rho_1-1)v
	\int_0^1 z^2\mu_1(\dd z)\int_0^1(1+zu)^{\rho_1-2}\dd u \cr
	\ar=:\ar
	\rho(1-\rho)\rho_1^2v^2 d_1
	-\rho\rho_1(\rho_1-1)v \bar{d}_1,\quad 0\le v\le1.
	\eeqlb
 Substituting  $z$ by $zv^{-1/\rho_1}$ we obtain
	\beqnn
	\ar\ar
	\int_1^\infty z^2\mu_1(\dd z)\int_0^1 H_1(v,zu)(1-u)\dd u \cr
	\ar=\ar
	v^{(\alpha_1-2)/\rho_1}\int_{v^{1/\rho_1}}^\infty z^2\mu_1(\dd z)\int_0^1 H_1(v,zv^{-1/\rho_1}u)(1-u)\dd u
	=:v^{\alpha_1/\rho_1}\tilde{H}_1(v)
	\eeqnn
	and
	\beqnn
	\ar\ar
	\int_1^\infty z^2\mu_1(\dd z)\int_0^1 H_2(v,zu)(1-u)\dd u \cr
	\ar=\ar
	v^{(\alpha_1-2)/\rho_1}\int_{v^{1/\rho_1}}^\infty z^2\mu_1(\dd z)\int_0^1 H_2(v,zuv^{-1/\rho_1})(1-u)\dd u
	=:v^{\alpha_1/\rho_1}\tilde{H}_2(v).
	\eeqnn
	It is elementary to see that
	\beqnn
	\lim_{v\to0}\tilde{H}_1(v)
	=\int_0^\infty z^2\mu_1(\dd z)\int_0^1
	\bar{H}_{1,\rho_1}(zu)(1-u)\dd u
	\eeqnn
	and
	\beqnn
	\lim_{v\to0}\tilde{H}_2(v)
	=\int_0^\infty z^2\mu_1(\dd z)\int_0^1
	\bar{H}_{2,\rho_1}(zu)(1-u)\dd u,
	\eeqnn
	where $\bar{H}_{1,\rho_1}(z):=z^{2\rho_1-2}(z^{\rho_1}+1)^{\rho-2}$
	and $\bar{H}_{2,\rho_1}(z):=z^{\rho_1-2}(z^{\rho_1}+1)^{\rho-1}$.
	Then
	\beqnn
	\int_0^1 z^2\mu_1(\dd z)\int_0^1
	\bar{H}_{2,\rho_1}(zu)(1-u)\dd u
	\le \int_0^1 z^{\rho_1}\mu_1(\dd z)
	=\frac{\alpha_1(\alpha_1-1)}
	{(\rho_1-\alpha_1)\Gamma(\alpha_1)\Gamma(2-\alpha_1)}
	\eeqnn
	and
	\beqnn
	\ar\ar
	\int_1^\infty z^2\mu_1(\dd z)\int_0^1
	\bar{H}_{2,\rho_1}(zu)(1-u)1_{\{zu\le1\}}\dd u \cr
	\ar\le\ar
	\int_1^\infty z^{\rho_1}\mu_1(\dd z)\int_0^{z^{-1}}u^{\rho_1-2}\dd u
	=(\rho_1-1)^{-1}\int_1^\infty z\mu_1(\dd z).
	\eeqnn
	By the dominated convergence, for $\rho\rho_1=\tilde{\rho}\in(0,1)$,
	\beqnn
	\ar\ar
	\lim_{\rho_1\to\infty}\int_1^\infty z^2\mu_1(\dd z)\int_0^1
	\bar{H}_{1,\rho_1}(zu)(1-u)1_{\{zu>1\}}\dd u\cr
	\ar=\ar
	\int_1^\infty z^2\mu_1(\dd z)\int_0^1
	(zu)^{\tilde{\rho}-2}(1-u)1_{\{zu>1\}}\dd u \cr
	\ar=\ar
	\int_0^1u^{\tilde{\rho}-2}(1-u)\dd u
	\int_{u^{-1}}^\infty z^{\tilde{\rho}}\mu_1(\dd z)
	=\frac{\alpha_1(\alpha_1-1)}
	{(\alpha_1-\tilde{\rho})\Gamma(\alpha_1)\Gamma(2-\alpha_1)}
	\int_0^1u^{\alpha_1-2}(1-u)\dd u
	\eeqnn
	and
	\beqnn
	\ar\ar
	\lim_{\rho_1\to\infty}\int_1^\infty z^2\mu_1(\dd z)\int_0^1
	\bar{H}_{2,\rho_1}(zu)(1-u)1_{\{zu>1\}}\dd u  \cr
	\ar=\ar
	\frac{\alpha_1(\alpha_1-1)}
	{(\alpha_1-\tilde{\rho})\Gamma(\alpha_1)\Gamma(2-\alpha_1)}
	\int_0^1u^{\alpha_1-2}(1-u)\dd u.
	\eeqnn
	Since $\rho(1-\rho)\rho_1^2-\rho\rho_1(\rho_1-1)
	=\rho\rho_1(1-\rho\rho_1)>0$,
	then there is a constant $\rho_0>0$ such that for all $\rho_1>\rho_0$,
	$\lim_{v\to0}\rho(1-\rho)\rho_1^2\tilde{H}_1(v)
	>\lim_{v\to0}\rho\rho_1(\rho_1-1)\tilde{H}_2(v)$.
	Now there is a constant $v_0:=v_0(\rho_1)\in(0,1)$ such that
	we have
	\beqlb\label{3.27}
	\rho(1-\rho)\rho_1^2\tilde{H}_1(v)\ge \rho\rho_1(\rho_1-1)\tilde{H}_2(v),
	\qquad 0<v\le v_0.
	\eeqlb
	For $v_0< v\le1$, we have
	$H_2(v,z)\le vv_0^{\rho-1}(1+z)^{\rho\rho_1-2}$ and then
	\beqnn
	\int_1^\infty z^2\mu_1(\dd z)\int_0^1 H_2(v,zu)(1-u)\dd u
	\le
	vv_0^{\rho-1} \int_0^1\dd u\int_1^\infty z^2(1+zu)^{\rho\rho_1-2}\mu_1(\dd z).
    \eeqnn
Replacing the variable $z$ with $y/u$ and using the fact
$\mu_1(\dd (y/u))=u^{\alpha_1}\mu_1(\dd y)$ we have
 \beqnn
 \ar\ar
\int_1^\infty z^2\mu_1(\dd z)\int_0^1 H_2(v,zu)(1-u)\dd u
\le
	vv_0^{\rho-1}\int_0^1u^{\alpha_1-2}\dd u\int_{ u}^\infty y^2(1+y)^{\rho\rho_1-2}\mu_1(\dd y) \cr
	\ar\ar\qquad\le
	vv_0^{\rho-1}\int_0^1u^{\alpha_1-2}\dd u
	\int_0^\infty z^2(1+z)^{\rho\rho_1-2}\mu_1(\dd z)=:
	vv_0^{\rho-1} c_0,\quad v_0<v\le1.
	\eeqnn
	Combining this with \eqref{3.8a} and \eqref{3.27} we obtain
	\beqnn
	\int_1^\infty K(v,z)\mu_1(\dd z)
	=
	\int_1^\infty z^2\mu_1(\dd z)\int_0^1
	K''_{zz}(v,zu) (1-u)\dd u
	\ge-\rho\rho_1(\rho_1-1)vv_0^{\rho-1} c_0
	\eeqnn
	for all $0<v\le1$.
	Together this with \eqref{3.8} one gets \eqref{3.28}.

	(ii)
	By \eqref{3.8a},
	\beqlb\label{3.23}
	K''_{zz}(v,z)
	\ar=\ar
	\rho\rho_1(1-\rho\rho_1)v^2(1+z)^{2\rho_1-2}\big[v[(1+z)^{\rho_1}-1]+1\big]^{\rho-2} \cr
	\ar\ar
	-\rho\rho_1(\rho_1-1)v(1-v)(1+z)^{\rho_1-2}\big[v[(1+z)^{\rho_1}-1]+1\big]^{\rho-2} \cr
	\ar=:\ar
	\rho\rho_1(1-\rho\rho_1)K_1(v,z)-\rho\rho_1(\rho_1-1)K_2(v,z).
	\eeqlb
	For $0\le v\le1$,
	\beqnn
	K_1(v,z)
	\ge
	v^2\big[(1+z)^{\rho_1}\big]^{\rho-2}
	(1+z)^{2\rho_1-2}
	=v^2(1+z)^{\rho\rho_1-2}.
	\eeqnn
	By \eqref{3.0a},
	\beqlb\label{3.24}
	\int_0^\infty z^2\mu_1(\dd z)\int_0^1K_1(v,uz)(1-u)\dd u
	\ge c_1(\rho\rho_1)v^2
	\eeqlb
	and for $\delta>0$
	\beqlb\label{3.26}
	\int_0^\delta z^2\mu_1(\dd z)\int_0^1K_2(v,uz)(1-u)\dd u
	\le
	v(1-v)\int_0^\delta z^2(1+z)^{\rho_1-1}\mu_1(\dd z)
	=:v(1-v)d_{1,\delta}.
	\eeqlb
	Moreover, for $0\le v\le 1$,
	\beqnn
	K_2(v,z)
	\le
	v(1-v)\big[v(1+z)^{\rho_1}\big]^{\rho-1}
	\big[(1-v)\big]^{-1}
	(1+z)^{\rho_1-1}
	=
	v^\rho
	(1+z)^{\rho\rho_1-2}
	\eeqnn
	and then for all $\delta>0$,
	\beqlb\label{3.25}
	\int_\delta^\infty z^2\mu_1(\dd z)\int_0^1K_2(v,uz)(1-u) \dd u
	\le
	v^\rho \int_\delta^\infty z^2 \mu_1(\dd z)\int_0^1(1+uz)^{\rho\rho_1-2}\dd u
	=:d_{2,\delta}v^\rho.
	\eeqlb
	By \eqref{6.0a} and \eqref{3.23},
	\beqnn
	K(v,z)
	\ar=\ar
	z^2\int_0^1K''_{zz}(v,zu)(1-u)\dd u \cr
	\ar=\ar
	\rho\rho_1(1-\rho\rho_1)z^2\int_0^1K_1(v,uz)(1-u)\dd u
	-\rho\rho_1(\rho_1-1)z^2\int_0^1K_2(v,uz)(1-u)\dd u.
	\eeqnn
	Combining this with \eqref{3.24}--\eqref{3.25} we obtain
	\eqref{3.11a}.
	\qed


	\subsubsection{Equations \eqref{1.1} with dominate drift terms}
\label{Equations {1.1} with dominate drifts}
	
In this subsubsection we consider equations \eqref{1.1}  under Condition \ref{c1} where the drifts dominate.
For $0<\rho<1$ and $\rho_1,\rho_2\ge1$ define
 \beqlb\label{3.7}
h(x,y)=-(x^{\rho_1}+y^{\rho_2})^\rho,\qquad x,y>0,
 \eeqlb
which is the key test function to prove Theorems \ref{t1.2}(i) and \ref{t1.4bb}
under Condition \ref{c1}.
In the following we estimate $\mathcal{L}h$ in Lemmas \ref{t7.2}
	and \ref{t7.3} for verifying assumption (iii) of Proposition \ref{t2.2}.
We first present a preliminary estimate on $\mathcal{L}h$ by \eqref{3.0} and Lemma \ref{t4.1a}. Recall \eqref{3.7bb}.
	By \eqref{3.10} and Lemma \ref{t4.1a}(i), for all large enough $\rho_1\ge 2$ and $\rho\rho_1\in(0,1)$, there is a constant $\tilde{d}_{1}>0$
	such that
	\beqnn
	\int_0^\infty K_z^1h(x,y)\mu_1(\dd z)
	\ar=\ar
	x^{-\alpha_1}(x^{\rho_1}+y^{\rho_2})^\rho
	\int_0^\infty K(\frac{x^{\rho_1}}{x^{\rho_1}+y^{\rho_2}},z)\mu_1(\dd z) \cr
	\ar\ge\ar
	-\rho\rho_1(\rho_1-1) \tilde{d}_{1}
	x^{\rho_1-\alpha_1}(x^{\rho_1}+y^{\rho_2})^{\rho-1}.
	\eeqnn
	Similarly, for all large enough $\rho_2\ge 2$ and $\rho\rho_2\in(0,1)$,
	there is a constant $\tilde{d}_{2}>0$ such that
	\beqnn
	\int_0^\infty K_z^2h(x,y)\mu_2(\dd z)
	\ge
	-\rho\rho_2(\rho_2-1) \tilde{d}_{2}
	y^{\rho_2-\alpha_2}(x^{\rho_1}+y^{\rho_2})^{\rho-1}.
	\eeqnn
 It thus follows 	from \eqref{3.0} that
	\beqlb\label{3.2}
	\ar\ar
	\mathcal{L}h(x,y) \cr
	\ar\ge\ar
	\rho(x^{\rho_1}+y^{\rho_2})^{\rho-1}
	\Big[b_{10}\rho_1x^{r_{10}-1+\rho_1}
	+b_{20}\rho_2y^{r_{20}-1+\rho_2}
	-a_1\rho_1x^{\theta_1-1+\rho_1}y^{\kappa_1}
	-a_2\rho_2 x^{\kappa_2}y^{\theta_2-1+\rho_2} \cr
	\ar\ar\qquad\qquad\qquad\quad
	-\rho_1(\rho_1-1)[b_{11}x^{r_{11}-2+\rho_1}
	+ b_{12}\tilde{d}_{1}x^{r_{12}-\alpha_1+\rho_1}] \cr
	\ar\ar\qquad\qquad\qquad\quad
	- \rho_2(\rho_2-1)[b_{21}y^{r_{21}-2+\rho_2}
	+ b_{22}\tilde{d}_{2}y^{r_{22}-\alpha_2+\rho_2}]\Big].
	\eeqlb
	By Lemma \ref{t4.1a}(ii),
	for $0<\rho\rho_1<1$ and $\rho_1>1$,
	\beqnn
	\ar\ar
	\int_0^\infty K_z^1h(x,y)\mu_1(\dd z)
	=x^{-\alpha_1}(x^{\rho_1}+y^{\rho_2})^\rho\int_0^\infty
	K(\frac{x^{\rho_1}}{x^{\rho_1}+y^{\rho_2}},z)\mu_1(\dd z)\cr
	\ar\ge\ar
	\rho\rho_1x^{-\alpha_1}\Big[(1-\rho\rho_1)c_1(\rho\rho_1)
	x^{2\rho_1}(x^{\rho_1}+y^{\rho_2})^{\rho-2} \cr
	\ar\ar\qquad\qquad
	-(\rho_1-1)d_{1,\delta}
	x^{\rho_1}y^{\rho_2}(x^{\rho_1}+y^{\rho_2})^{\rho-2}
	-(\rho_1-1)d_{2,\delta}
	x^{\rho\rho_1} \Big].
	\eeqnn
	Similar results hold for $\int_0^\infty K_z^2h(x,y)\mu_2(\dd z)$ by symmetry.
	Thus, we also have
	\beqlb\label{3.2bb}
	\ar\ar
	\mathcal{L}h(x,y) \cr
	\ar\ge\ar
	\rho(x^{\rho_1}+y^{\rho_2})^{\rho-1}
	\Big[b_{10}\rho_1x^{r_{10}-1+\rho_1}
	-a_1\rho_1x^{\theta_1-1+\rho_1}y^{\kappa_1}
	-a_2\rho_2 x^{\kappa_2}y^{\theta_2-1+\rho_2} \cr
	\ar\ar\qquad\qquad\qquad\quad
	-\rho_1(\rho_1-1)[b_{11}x^{r_{11}-2+\rho_1}
	+ b_{12}\tilde{d}_{1}x^{r_{12}-\alpha_1+\rho_1}] \Big] \cr
	\ar\ar
	+\rho\rho_2(x^{\rho_1}+y^{\rho_2})^{\rho-2}\Big[(1-\rho\rho_2)[b_{21}y^{r_{21}-2+2\rho_2}
	+b_{22}c_2(\rho\rho_2)y^{r_{22}-\alpha_2+2\rho_2}] \cr
	\ar\ar\qquad\qquad\qquad\qquad\quad
	-(\rho_2-1)x^{\rho_1}[b_{21}y^{r_{21}-2+\rho_2}
	+b_{22}\tilde{d}_{1,\delta}y^{r_{22}-\alpha_2+\rho_2}]\Big] \cr
	\ar\ar
	-\rho\rho_2(\rho_2-1)b_{22}\tilde{d}_{2,\delta}
	y^{r_{22}-\alpha_2+\rho\rho_2},
	\eeqlb
	where $\tilde{d}_{1,\delta},\tilde{d}_{2,\delta}>0$
satisfy $\lim_{\delta\to\infty}\tilde{d}_{2,\delta}=0$.

Based on the above estimation we state the following lemma, whose
proof utilizes Lemmas \ref{t2.3} and \ref{t7.1}.

	\blemma\label{t7.2}
	Suppose that $\theta_1-1<r_1<0$, $\theta_2-1<r_2<0$
	and the assumptions of Theorem \ref{t1.2}(i) hold. Then
	there are constants $\rho_1,\rho_2,\rho>0$ and $C,c>0$ such that
	\beqlb\label{7.1}
	\mathcal{L}h(x,y)\ge C,
	\qquad 0<x,y\le c.
	\eeqlb
	\elemma
	\proof
	We first establish the proof under Condition \ref{c1}(i).
	Let $\rho_1,\rho_2\ge2$ be large enough constants satisfying \eqref{3.6}.
	Let $0<\rho<1$ be small enough such that
	$r_{10}-1+\rho_1\rho\le0$ and $r_{20}-1+\rho_2\rho\le0$.
	By \eqref{3.2},
	\beqlb\label{3.9}
	\mathcal{L}h(x,y)
	\ge
	\rho(x^{\rho_1}+y^{\rho_2})^{\rho-1}
	[J_1(x,y)-J_2(x,y)-J_3(x)-J_4(y)],\quad 0<x,y<1,
	\eeqlb
	where
	\beqlb\label{3.4b}
	\begin{aligned}	
	J_1(x,y)&:=
	b_{10}\rho_1x^{r_{10}-1+\rho_1}
	+b_{20}\rho_2y^{r_{20}-1+\rho_2},\\
	J_2(x,y)&:= a_1\rho_1x^{\theta_1-1+\rho_1}y^{\kappa_1}
	+a_2\rho_2 x^{\kappa_2}y^{\theta_2-1+\rho_2}
   \end{aligned}
	\eeqlb
	and
	\begin{equation}\label{3.4c}
		\begin{aligned}
			J_3(x)
			&:=
			\rho_1(\rho_1-1)[b_{11}x^{r_{11}-2+\rho_1}
			+ b_{12}\tilde{d}_1x^{r_{12}-\alpha_1+\rho_1}],\cr
			J_4(y)
			&:=
			\rho_2(\rho_2-1)[b_{21}y^{r_{21}-2+\rho_2}
			+ b_{22}\tilde{d}_2y^{r_{22}-\alpha_2+\rho_2}]
		\end{aligned}
	\end{equation}
	for some constants $\tilde{d}_1,\tilde{d}_2>0$.
	By Lemma \ref{t7.1},
	there is a constant $0<c_1<1$ such that
	\beqlb\label{3.4aa}
	4^{-1}J_1(x,y)\ge J_2(x,y),\qquad
	0<x,y\le c_1.
	\eeqlb
	By Condition \ref{c1}(i), there is a constant $0<c_2\le c_1$ such that
	$8^{-1}J_1(x,y)\ge J_3(x)$ and $8^{-1}J_1(x,y)\ge J_4(y)$
	for all $0<x,y\le c_2$.
	Combining these with \eqref{3.9} and \eqref{3.4aa} we obtain
	$$\mathcal{L}h(x,y)\ge2^{-1}\rho(x^{\rho_1}+y^{\rho_2})^{\rho-1}J_1(x,y)$$
	for all $0<x,y\le c_2$.
	By Lemma \ref{t2.3}(ii), we have
	\beqlb\label{3.44b}
	J_1(x,y)
	\ar=\ar
	b_{10}\rho_1x^{r_{10}-1+\rho_1\rho}x^{\rho_1(1-\rho)}
	+b_{20}\rho_2y^{r_{20}-1+\rho_2\rho}y^{\rho_2(1-\rho)} \cr
	\ar\ge\ar
	C
	[x^{\rho_1(1-\rho)}+y^{\rho_2(1-\rho)}]
	\ge C(x^{\rho_1}+y^{\rho_2})^{1-\rho}
	\eeqlb
with $C:=(b_{10}\rho_1)\wedge(b_{20}\rho_2)$.
	Thus,
$\mathcal{L}h(x,y)
	\ge
	2^{-1}\rho C$
for all $0<x,y\le c_2$,
	which gives \eqref{7.1}.

	Since the proofs under Condition \ref{c1}(ii)-(iii) are similar,
	we only present that under Condition \ref{c1}(ii).
	If Condition \ref{c1}(ii) and inequality \eqref{1.3} hold,
	there are large enough constants $\rho_1,\rho_2\ge2$ such that \eqref{3.6} holds
	and
	$\rho_1/\rho_2<r_1/r_2$.
	Moreover, for all $0<\rho<1$,
	\beqlb\label{3.23b}
	\frac{\rho_1}{r_1+2\rho_1}
	+\frac{r_2+\rho_2}{r_2+2\rho_2}>1,\quad
	\frac{\rho_1(2-\rho)}{r_1+2\rho_1}
	+\frac{r_2+\rho\rho_2}{r_2+2\rho_2}>1.
	\eeqlb
	In the following, let $0<\rho<1$ be small enough such that
	$r_{10}-1+\rho_1\rho=r_1+\rho_1\rho\le0$ and $r_2+\rho_2\rho\le0$.
	Recall \eqref{3.4b} and \eqref{3.4c}.
	By \eqref{3.2bb},
	\beqnn
	\mathcal{L}h(x,y)
	\ar\ge\ar
	\rho(x^{\rho_1}+y^{\rho_2})^{\rho-1}
	\big[b_{10}\rho_1x^{r_1+\rho_1}-J_2(x,y)
	-J_3(x) \big] \cr
	\ar\ar
	+\rho(x^{\rho_1}+y^{\rho_2})^{\rho-2}\big[\tilde{c}_1y^{r_2+2\rho_2}
	-\tilde{c}_2x^{\rho_1}y^{r_2+\rho_2}\big]
	-\tilde{c}_0\tilde{d}_{2,\delta}
	y^{r_2+\rho\rho_2},
	\eeqnn
	where $\tilde{c}_0:=\rho\rho_2(\rho_2-1)b_{22}$,
	$\tilde{c}_1:=\rho_2(1-\rho\rho_2)
	[b_{21}\wedge (b_{22}c_2(\rho\rho_2))]$
	and
	$\tilde{c}_2:=\rho_2(\rho_2-1)[b_{21}
	+b_{22}\tilde{d}_{1,\delta}]$.
	Then for $0<\sigma_1<b_{10}\rho_1$ and $0<\sigma_2<\tilde{c}_1$,
	\beqlb\label{3.11}
	\mathcal{L}h(x,y)
	\ar\ge\ar
	\rho(x^{\rho_1}+y^{\rho_2})^{\rho-1}
	\big[I_1(x,y)-J_2(x,y)
	-J_3(x) \big] \cr
	\ar\ar
	+\rho(x^{\rho_1}+y^{\rho_2})^{\rho-2}\big[
	I_2(x,y)
	-(\tilde{c}_2+\sigma_2)x^{\rho_1}y^{r_2+\rho_2}\big]
	-\tilde{c}_0\tilde{d}_{2,\delta}
	y^{r_2+\rho\rho_2},
	\eeqlb
	where $I_1(x,y):=(b_{10}\rho_1-\sigma_1)x^{r_1+\rho_1}+\sigma_2 y^{r_2+\rho_2}$
	and $I_2(x,y):=\sigma_1x^{r_1+2\rho_1}
	+(\tilde{c}_1-\sigma_2)y^{r_2+2\rho_2}$.
	Since $\lim_{\delta\to\infty}\tilde{d}_{2,\delta}=0$,
 there is a constant $\delta>0$ such that $I_2(x,y)\ge
	8\tilde{c}_0\tilde{d}_{2,\delta}y^{r_2+2\rho_2}$.
	By Lemma \ref{t7.1}, \eqref{3.6} and \eqref{3.23b},
	there is a constant $0<c_3<1$ such that
	for all $0<x,y\le c_3$ we have
	\beqlb\label{3.4aaa}
	I_1(x,y)\ge 4J_2(x,y),~
	I_2(x,y)
	\ge2(\tilde{c}_2+\sigma_2)x^{\rho_1}y^{r_2+\rho_2},~
	I_2(x,y)
	\ge8\tilde{c}_0\tilde{d}_{2,\delta}x^{\rho_1(2-\rho)}y^{r_2+\rho\rho_2}.
	\eeqlb
	It follows from Lemma \ref{t2.3}(ii) that
	\beqlb\label{3.4ab}
	I_2(x,y)
	-2\tilde{c}_0\tilde{d}_{2,\delta}
	y^{r_2+\rho\rho_2}(x^{\rho_1}+y^{\rho_2})^{2-\rho}
	\ge
	I_2(x,y)
	-4\tilde{c}_0\tilde{d}_{2,\delta}
	[y^{r_2+\rho\rho_2}x^{\rho_1(2-\rho)}
	+y^{r_2+2\rho_2}] \ge0
	\eeqlb
	for all $0<x,y\le c_3$.
	By \eqref{1.8}, there is a constant $0<c_4<c_3$
such that
	$4^{-1}I_1(x,y)\ge J_3(x)$ for all $0<x,y\le c_4$.
	Combining this with \eqref{3.11}--\eqref{3.4ab}, one obtains
	\beqnn
	\mathcal{L}h(x,y)
	\ge
	2^{-1}\rho(x^{\rho_1}+y^{\rho_2})^{\rho-1}
	I_1(x,y),\quad
	0<x,y\le c_4.
	\eeqnn
Similarly to \eqref{3.44b}, we have
	$I_1(x,y)\ge [(b_{10}\rho_1-\sigma_1)\wedge\sigma_2]
	(x^{\rho_1}+y^{\rho_2})^{1-\rho}$,
	which gives \eqref{7.1}.
	\qed

	Let $\delta_0>0$ be the constant
	determined in Lemma \ref{t5.2}.
	For $\rho_1,\rho_2,\rho>0$ define
	\beqlb\label{3.7b}
	\tilde{h}(x,y):=-(\delta_0 x^{\rho_1}+y^{\rho_2})^\rho,\qquad x,y>0.
	\eeqlb
Similarly to Lemma \ref{t7.2}, we have the following estimate
of $\mathcal{L}\tilde{h}$.
	\blemma\label{t7.3}
According to  the assumptions of Theorem \ref{t1.4bb},
	there are constants $\rho_1,\rho_2,\delta_0,\rho>0$ and $C,c>0$ such that
	\beqlb\label{7.1b}
	\mathcal{L}\tilde{h}(x,y)\ge C,
	\qquad 0<x,y\le c.
	\eeqlb
	\elemma
	\proof
	Using Lemma \ref{t5.2},
	the proof is a modification of that of Lemma \ref{t7.2},
and we omit it here.
	\qed

	\subsubsection{Equations \eqref{1.1} with dominant diffusion or jump terms}\label{Equations {1.1} with dominate diffusive or jump terms}
	
 In this subsubsection we study the SDE system \eqref{1.1} with dominant diffusion or jump terms, that is, under the assumption \eqref{1.6}.
The key function $h$ is defined in Lemma \ref{t7.4} under Condition \ref{c2}(i),
in Lemma \ref{t7.5} under Condition \ref{c2}(ii),
and in Lemma \ref{t7.7} under assumption (iii) of Theorem \ref{t1.2}.
The estimates of $\mathcal{L}h$,  key to verifying assumption (iii) of Proposition \ref{t2.2}, are given in Lemmas \ref{t7.4}, \ref{t7.5} and \ref{t7.7} for different cases. Applying these estimates and Proposition \ref{t2.2}, the proofs of Theorems \ref{t1.2} and \ref{t1.4bb} are given at the end of this subsubsection.
	For $\varepsilon\in(0,1)$ define a nonnegative bounded function
	\beqlb\label{7.5}
	\tilde{h}(y):=y-(1+\varepsilon)^{-1}y^{1+\varepsilon}(1+y)^{-\varepsilon},\qquad y>0.
	\eeqlb
	Then
	\beqlb\label{3.13}
	\tilde{h}'(y)=1-y^{\varepsilon}(1+y)^{-\varepsilon}
	+\varepsilon(1+\varepsilon)^{-1}y^{1+\varepsilon}(1+y)^{-1-\varepsilon}\ge0,~
	\tilde{h}''(y)
	=
	-\varepsilon y^{\varepsilon-1}(1+y)^{-2-\varepsilon}.
	\eeqlb
	Thus, there is a constant $c_0=c_0(\varepsilon)\in(0,1)$ such that
	\beqlb\label{3.13b}
	\tilde{h}(y)\ge 2^{-1}y,~2^{-1}\le\tilde{h}'(y)\le1,~
	-\tilde{h}''(y)\ge 2^{-1}\varepsilon y^{\varepsilon-1},
	\qquad 0<y\le2c_0.
	\eeqlb
	It follows that for all $0<y\le c_0$ we have
	\beqlb\label{3.29}
	\tilde{h}(y+yzu)\le y+yzu\le 2y\le4\tilde{h}(y),
	\qquad 0<y\le c_0,~0<z,u\le 1.
	\eeqlb
	\blemma\label{t7.4}
	Suppose that \eqref{1.6} holds for $i=1,2$,
	$\theta_1-1<r_1<0$,
	$\theta_2-1<r_2<0$ and \eqref{1.3} and that condition \ref{c2}(i) holds.
	In addition, under \eqref{3.113} let
	\beqlb\label{7.3}
	h(x,y):=-[x^{\rho_1}+\tilde{h}(y)]^\rho,\qquad x,y>0,
	\eeqlb
	and under \eqref{3.113b} let
	\beqnn
	h(x,y):=-[y^{\rho_2}+\tilde{h}(x)]^\rho,\qquad x,y>0.
	\eeqnn
	Then there are constants $\rho_1,\rho_2,\rho,\varepsilon>0$
	and $C,c>0$ such that $\mathcal{L}h(x,y)\ge C$ for all
	$0<x,y\le c$.
	\elemma

Before proving Lemma \ref{t7.4}, we first show a preliminary estimate
of $\mathcal{L}h$ by \eqref{3.0} in the following lemma.
\blemma\label{t7.4b}
Suppose that \eqref{1.6} holds for $i=1,2$.
Let $h$ be the function defined in \eqref{7.3}
for $0<\rho,\varepsilon<1$ and $\rho_1>1$.
Then there are constants $\tilde{c}_1,\tilde{c}_2,\tilde{c}_3,\tilde{c}_4>0$ such that
	\beqlb\label{3.16}
	\mathcal{L}h(x,y)
	\ar\ge\ar
	-8^{-1}\tilde{c}_2\rho x^{r_1+\rho\rho_1}
	+\rho[x^{\rho_1}+\tilde{h}(y)]^{\rho-1}
	\big[\tilde{c}_1 y^{r_2+1+\varepsilon}
	-a_1\rho_1x^{\rho_1-1+\theta_1}y^{\kappa_1}
	-2a_2x^{\kappa_2}y^{\theta_2}\big] \cr
	\ar\ar
	+\rho[x^{\rho_1}+\tilde{h}(y)]^{\rho-2}
	\big[\tilde{c}_2 x^{r_1+2\rho_1}+\tilde{c}_3y^{r_2+2}
	-\tilde{c}_4x^{r_1+\rho_1}\tilde{h}(y)\big]
	,\quad
	0<x,y\le c_0,
	\eeqlb
where $c_0\in(0,1)$ is the constant appearing in \eqref{3.13b}.
\elemma
\proof
	Observe that
	\beqlb\label{3.12a}
	h_x'(x,y)=-\rho\rho_1[x^{\rho_1}+\tilde{h}(y)]^{\rho-1}x^{\rho_1-1},\quad
	h_y'(x,y)=-\rho[x^{\rho_1}+\tilde{h}(y)]^{\rho-1}
	\tilde{h}'(y)
	\eeqlb
	and
	\beqlb\label{3.12b}
	h_{xx}''(x,y)
	=
	\rho\rho_1[x^{\rho_1}+\tilde{h}(y)]^{\rho-2}
	[(1-\rho\rho_1)x^{2\rho_1-2}-(\rho_1-1)x^{\rho_1-2}\tilde{h}(y)].
	\eeqlb
	Moreover, by \eqref{3.13b},
	\beqlb\label{3.12}
	h_{yy}''(x,y)
	\ar=\ar
	\rho(1-\rho)[x^{\rho_1}+\tilde{h}(y)]^{\rho-2}
	|\tilde{h}'(y)|^2
	-\rho [x^{\rho_1}+\tilde{h}(y)]^{\rho-1} \tilde{h}''(y) \cr
	\ar\ge\ar
	2^{-2}\rho(1-\rho)[x^{\rho_1}+\tilde{h}(y)]^{\rho-2}
	+2^{-1}\varepsilon\rho y^{\varepsilon-1}[x^{\rho_1}+\tilde{h}(y)]^{\rho-1},
	~~0<y\le 2c_0.
	\eeqlb
	Let $K(v,z)$ be the function defined in \eqref{3.7bb}.
	By \eqref{3.10} and \eqref{3.11a}, we have
	\beqlb\label{3.14}
	\ar\ar
	\int_0^\infty K_z^1h(x,y)\mu_1(\dd z)
	=
	x^{-\alpha_1}[x^{\rho_1}+\tilde{h}(y)]^\rho
	\int_0^\infty K(\frac{x^{\rho_1}}{x^{\rho_1}+\tilde{h}(y)},z)\mu_1(\dd z) \cr
	\ar\ar\quad\ge
	\rho\rho_1x^{-\alpha_1}\Big[(1-\rho\rho_1)c_1(\rho\rho_1)
	x^{2\rho_1}[x^{\rho_1}+\tilde{h}(y)]^{\rho-2} \cr
	\ar\ar\qquad\qquad\qquad
	-(\rho_1-1)d_{1,\delta}
	x^{\rho_1}\tilde{h}(y)[x^{\rho_1}+\tilde{h}(y)]^{\rho-2}
	-(\rho_1-1)d_{2,\delta}
	x^{\rho\rho_1} \Big],
	\eeqlb
	where the constants $d_{1,\delta},d_{2,\delta}>0$
	satisfy
	$\lim_{\delta\to\infty}d_{2,\delta}=0$.

	Combining with \eqref{3.13b}, \eqref{3.29} and \eqref{3.12},
	for $0<z,u\le1$ and $0<y\le c_0$ we have
	\beqlb\label{3.12b}
	h_{yy}''(x,y+yzu)
	\ar\ge\ar
	2^{-2}\rho(1-\rho)[x^{\rho_1}+4\tilde{h}(y)]^{\rho-2}
	+2^{-2}\rho\varepsilon  [x^{\rho_1}+4\tilde{h}(y)]^{\rho-1} y^{\varepsilon-1} \cr
	\ar\ge\ar
	2^{2\rho-6}\rho(1-\rho)[x^{\rho_1}+\tilde{h}(y)]^{\rho-2}
	+2^{2\rho-4}\rho\varepsilon [x^{\rho_1}+\tilde{h}(y)]^{\rho-1} y^{\varepsilon-1}.
	\eeqlb
	Since $h_{yy}''(x,y)\ge0$ for all $y>0$,
	then by \eqref{6.0a},
	\beqnn
	\int_0^{\infty}K_{yz}^2h(x,y)\mu_2(\dd z)
	\ar=\ar
	y^2\int_0^{\infty}z^2 \mu_2(\dd z)\int_0^1h_{yy}''(x,y+yzu)(1-u)\dd u \cr
	\ar\ge\ar
	y^2\int_0^1z^2 \mu_2(\dd z)\int_0^1h_{yy}''(x,y+yzu)(1-u)\dd u \cr
	\ar\ge\ar
	\hat{c}_1\rho\varepsilon y^2 \big[(1-\rho)[x^{\rho_1}+\tilde{h}(y)]^{\rho-2}
	+ [x^{\rho_1}+\tilde{h}(y)]^{\rho-1} y^{\varepsilon-1}\big],
	~~0<y\le c_0,
	\eeqnn
for $\hat{c}_1:=2^{2\rho-6}\int_0^1z^2\mu_2(\dd z)$,
where \eqref{3.12b} is used in the last inequality.
Then by \eqref{6.0b},
	\beqlb\label{3.15}
	\ar\ar
	\int_0^\infty K_{z}^2h(x,y)\mu_2(\dd z)
	=y^{-\alpha_2}\int_0^\infty K_{y z}^2h(x,y)\mu_2(\dd z)
	\cr
	\ar\ar\qquad\ge
	\hat{c}_1\rho\varepsilon y^{2-\alpha_2} \big[(1-\rho)[x^{\rho_1}+\tilde{h}(y)]^{\rho-2}
	+ [x^{\rho_1}+\tilde{h}(y)]^{\rho-1} y^{\varepsilon-1}\big],
	~~0<y\le c_0.
	\eeqlb
	Further combining \eqref{3.0}
	with \eqref{3.12a}--\eqref{3.14} and \eqref{3.15} one obtains
	\beqnn
	\ar\ar
	\mathcal{L}h(x,y) \cr
	\ar\ge\ar
	\rho[x^{\rho_1}+\tilde{h}(y)]^{\rho-1}
	\big[ b_{21} 2^{-1}\varepsilon y^{r_{21}-1+\varepsilon}
	+b_{22} \hat{c}_1 \varepsilon y^{r_{22}+1-\alpha_2+\varepsilon}
	-a_1\rho_1x^{\rho_1-1+\theta_1}y^{\kappa_1}
	-a_2\tilde{h}'(y)x^{\kappa_2}y^{\theta_2}\big] \cr
	\ar\ar
	+\rho [x^{\rho_1}+\tilde{h}(y)]^{\rho-2}\Big[
	b_{11}\rho_1(1-\rho\rho_1)x^{r_{11}+2\rho_1-2}
	+ b_{12}\rho_1(1-\rho\rho_1)c_1(\rho\rho_1)
	x^{r_{12}+2\rho_1-\alpha_1} \cr
	\ar\ar\qquad\qquad\qquad\qquad~
	+b_{21}(1-\rho)
	2^{-2}y^{r_{21}}
	+b_{22}(1-\rho)\hat{c}_1 \varepsilon y^{r_{22}+2-\alpha_2}
	\cr
	\ar\ar\qquad\qquad\qquad\qquad~
	-b_{11}\rho_1(\rho_1-1) x^{r_{11}+\rho_1-2}\tilde{h}(y)
	-b_{12}\rho_1(\rho_1-1) d_{1,\delta}
	x^{r_{12}+\rho_1-\alpha_1}\tilde{h}(y)\Big] \cr
	\ar\ar
	-b_{12}\rho\rho_1 (\rho_1-1)d_{2,\delta}
	x^{r_{12}+\rho\rho_1-\alpha_1},\qquad
	0<x,y\le c_0.
	\eeqnn
	Since $\lim_{\delta\to\infty}d_{2,\delta}=0$,
	there is a large enough constant $\delta>0$ such that
	$b_{12}\rho_1 (\rho_1-1)d_{2,\delta}\le 8^{-1}\tilde{c}_2$
	with $\tilde{c}_2:=\rho_1(1-\rho\rho_1) [b_{11}\wedge
	(b_{12}c_1(\rho\rho_1))]$.
	By \eqref{3.13}, $\tilde{h}'(y)\le 2$ for all $y>0$ we get
\eqref{3.16} with $\tilde{c}_1:=\varepsilon [(b_{21} 2^{-1})
	\wedge (b_{22} \hat{c}_1)],\tilde{c}_3:=[b_{21}(1-\rho)
	2^{-2}]
	\wedge[b_{22}\hat{c}_1 \varepsilon]$
	and
	$\tilde{c}_{4}:=\rho_1(\rho_1-1)[b_{11}
	+b_{12} d_{1,\delta}]$.
This ends the proof.
\qed

Now we are ready to complete the proof of Lemma \ref{t7.4}
applying Lemma \ref{t2.3}(ii)--(iii).

\noindent{\it Proof of the Lemma \ref{t7.4}.}
	Since the proofs are similar, we only present that under \eqref{3.113}.
	By \eqref{3.113} and \eqref{1.3}, we obtain
	\beqnn
	1<\frac{\kappa_2(r_2+1)}{r_2+1-\theta_2}-r_1,\quad \frac{(r_1+1-\theta_1)}{\kappa_1}<\frac{\kappa_2}{r_2+1-\theta_2},\quad
	\frac{r_1}{r_2}<\frac{\kappa_2(r_2+1)}{r_2+1-\theta_2}-r_1.
	\eeqnn
	Then there is a constant $\rho_1$ such that
	\beqnn
	1<\rho_1,\quad \frac{r_1}{r_2}<\rho_1,\quad
	\frac{(r_1+1-\theta_1)(r_2+1)}{\kappa_1}-r_1<\rho_1
	<\frac{\kappa_2(r_2+1)}{r_2+1-\theta_2}-r_1,
	\eeqnn
	which gives
	\beqnn
	\frac{\rho_1-1+\theta_1}{r_1+\rho_1}+\frac{\kappa_1}{r_2+1}>1,\qquad
	\frac{\kappa_2}{r_1+\rho_1}+\frac{\theta_2}{r_2+1}>1
	\eeqnn
	and
	\beqlb\label{3.111}
	\frac{r_1+\rho_1}{r_1+2\rho_1}+\frac{1}{r_2+2}>1,\qquad
	\frac{r_1+\rho\rho_1}{r_1+2\rho_1}+\frac{2-\rho}{r_2+2}>1
	\eeqlb
	for small enough $0<\rho<1$ with
	\beqlb\label{3.17}
	r_1+\rho\rho_1<0,\qquad r_2+\rho<0.
	\eeqlb
	Moreover, there is a sufficiently small constant $0<\varepsilon<1$ such that
	\beqlb\label{3.112}
	\frac{\rho_1-1+\theta_1}{r_1+\rho_1}+\frac{\kappa_1}{r_2+1+\varepsilon}>1,\qquad
	\frac{\kappa_2}{r_1+\rho_1}+\frac{\theta_2}{r_2+1+\varepsilon}>1.
	\eeqlb

	By Lemma \ref{t2.3}(ii),
	\beqnn
	x^{r_1+\rho\rho_1}
	=[x^{\rho_1}+\tilde{h}(y)]^{\rho-2}
	x^{r_1+\rho\rho_1}[x^{\rho_1}+\tilde{h}(y)]^{2-\rho}
	\le
	2[x^{\rho_1}+\tilde{h}(y)]^{\rho-2}
	[x^{r_1+2\rho_1}+x^{r_1+\rho\rho_1}\tilde{h}(y)^{2-\rho}]
	\eeqnn
	and
	\beqnn
	x^{r_1+\rho_1}
	=[x^{\rho_1}+\tilde{h}(y)]^{-1}[x^{r_1+2\rho_1}+x^{r_1+\rho_1}\tilde{h}(y)].
	\eeqnn
 Thus, it follows from Lemma \ref{t7.4b} that
	\beqlb\label{3.114}
	\mathcal{L}h(x,y)
	\ar\ge\ar
	\rho[x^{\rho_1}+\tilde{h}(y)]^{\rho-1}
	\big[4^{-1}\tilde{c}_2x^{r_1+\rho_1}+\tilde{c}_1 y^{r_2+1+\varepsilon}
	-a_1\rho_1x^{\rho_1-1+\theta_1}y^{\kappa_1}
	-2a_2x^{\kappa_2}y^{\theta_2}\big] \cr
	\ar\ar
	+\rho [x^{\rho_1}+\tilde{h}(y)]^{\rho-2}
	\big[2^{-1}\tilde{c}_2 x^{r_1+2\rho_1}+\tilde{c}_3y^{r_2+2}
	-(\tilde{c}_4+4^{-1}\tilde{c}_2)x^{r_1+\rho_1}\tilde{h}(y) \cr
	\ar\ar\qquad\qquad\qquad\qquad
	-4^{-1}\tilde{c}_2x^{r_1+\rho\rho_1}\tilde{h}(y)^{2-\rho}\big] \cr
	\ar\ge\ar
	\rho[x^{\rho_1}+\tilde{h}(y)]^{\rho-1}
	\big[4^{-1}\tilde{c}_2x^{r_1+\rho_1}+\tilde{c}_1 y^{r_2+1+\varepsilon}
	-a_1\rho_1x^{\rho_1-1+\theta_1}y^{\kappa_1}
	-2a_2x^{\kappa_2}y^{\theta_2}\big] \cr
	\ar\ar
	+\rho [x^{\rho_1}+\tilde{h}(y)]^{\rho-2}
	\big[2^{-1}\tilde{c}_2 x^{r_1+2\rho_1}+\tilde{c}_3y^{r_2+2}
	-(\tilde{c}_4+4^{-1}\tilde{c}_2)x^{r_1+\rho_1}y \cr
	\ar\ar\qquad\qquad\qquad\qquad
	-4^{-1}\tilde{c}_2x^{r_1+\rho\rho_1}y^{2-\rho}\big] \cr
	\ar=:\ar
	\rho[x^{\rho_1}+\tilde{h}(y)]^{\rho-1} I_1(x,y)
	+\rho[x^{\rho_1}+\tilde{h}(y)]^{\rho-2}I_2(x,y),\quad
	0<x,y\le c_0,
	\eeqlb
	where we use the fact $\tilde{h}(y)\le y$ for all $y>0$
in the second inequality.
By Lemma \ref{t2.3}(iii), \eqref{3.111} and \eqref{3.112}, there is a constant
	$0<c_1\le c_0$ such that
	\beqnn
	I_1(x,y)\ge0,~~
	I_2(x,y)\ge
	4^{-1}[\tilde{c}_2 x^{r_1+2\rho_1}+\tilde{c}_3y^{r_2+2}],
	~\qquad 0<x,y<c_1.
	\eeqnn
 Thus, it follows from \eqref{3.17}, \eqref{3.114},
	the fact $\tilde{h}(y)\le y$, and Lemma \ref{t2.3}(ii) that
	\beqnn
	\mathcal{L}h(x,y)
	\ar\ge\ar
	4^{-1}\rho[x^{\rho_1}+y]^{\rho-2}
	[\tilde{c}_2 x^{\rho_1(2-\rho)} x^{r_1+\rho\rho_1}
	+\tilde{c}_3y^{2-\rho}y^{r_2+\rho}] \cr
	\ar\ge\ar
	4^{-1}\rho(\tilde{c}_2 \wedge\tilde{c}_3)[x^{\rho_1}+y]^{\rho-2}
	[x^{\rho_1(2-\rho)}+y^{2-\rho}]
	\ge8^{-1}\rho(\tilde{c}_2 \wedge\tilde{c}_3),
	\quad 0<x,y<c_1,
	\eeqnn
	which completes the proof.
	\qed

	\blemma\label{t7.5}
 Given $\theta_1-1<r_1<0$ and $\theta_2-1<r_2<0$, under  Condition \ref{c2}(ii), suppose that  \eqref{1.3} holds and \eqref{1.6} holds for $i=1,2$.
	Under \eqref{3.122a} let
	\beqlb\label{7.4}
	\hat{h}(x,y):=(x+y^\delta)^{\rho_1}+\tilde{h}(y),~~
	h(x,y):=-\hat{h}(x,y)^\rho,\qquad x,y>0,
	\eeqlb
	and under \eqref{3.122b} let
	\beqnn
	\hat{h}(x,y):=(y+x^\delta)^{\rho_2}+\tilde{h}(x),~~
	h(x,y):=-\hat{h}(x,y)^\rho, \qquad x,y>0.
	\eeqnn
	Then there are constants $\rho_1,\rho_2,\rho,\varepsilon>0$
	and $C, c>0$ such that $\mathcal{L}h(x,y)\ge C$ for all
	$0<x,y\le c$.
	\elemma

Before proving Lemma \ref{t7.5}, we first give an estimate
of $\mathcal{L}h$ by \eqref{3.0} in the following lemma.
\blemma\label{t7.5b}
Let $0<\rho,\rho_1,\varepsilon<1$ and $\delta>1$ satisfy
$\rho_1\delta>1$, $0<\varepsilon<\delta\rho_1-1$ and
$\delta\rho\rho_1<1$.
Let $h$ be the function defined by \eqref{7.4}.
Then there are constants $c_1,\tilde{c}_1,\tilde{c}_2,\tilde{c}_3>0$
such that
	\beqnn
	\mathcal{L}h(x,y)
\ge
	\rho \hat{h}(x,y)^{\rho-1}[I_1(x,y)-I_2(x,y)]
	-b_{22}\tilde{c}_1y^{r_{22}},\qquad 0<x,y\le c_1,
	\eeqnn
	where
	\beqnn
	I_1(x,y):=
	\tilde{c}_2\big[(x+y^\delta)^{\rho_1-2}x^{r_{11}}
	+(x+y^\delta)^{\rho_1-\alpha_1}x^{r_{12}}\big]
	+\tilde{c}_3 y^{\varepsilon+1}[y^{r_{21}-2}+y^{r_{22}-\alpha_2}]
	\eeqnn
	and
	\beqnn
	I_2(x,y):=a_1\rho_1x^{\theta_1}y^{\kappa_1}(x+y^\delta)^{\rho_1-1}
	+2a_2 y^{\theta_2}x^{\kappa_2}.
	\eeqnn
\elemma
\proof
	Observe that
	\beqlb\label{3.20}
	h_x'(x,y)=-\rho\rho_1\hat{h}(x,y)^{\rho-1}(x+y^\delta)^{\rho_1-1},~
	h_{xx}''(x,y)
	\ge
	\rho\rho_1(1-\rho_1)\hat{h}(x,y)^{\rho-1}(x+y^\delta)^{\rho_1-2}
	\eeqlb
	and
	\beqnn
	h_y'(x,y)
	\ar=\ar
	-\rho \hat{h}(x,y)^{\rho-1}
	\big[\rho_1\delta(x+y^\delta)^{\rho_1-1}y^{\delta-1}+\tilde{h}'(y)\big], \cr
	h_{yy}''(x,y)
	\ar\ge\ar
	\rho \hat{h}(x,y)^{\rho-1}
	\big[
	-\rho_1(\delta-1)\delta(x+y^\delta)^{\rho_1-1}y^{\delta-2}
	-\tilde{h}''(y) \big].
	\eeqnn
	Since $\rho_1\delta>1$ and $0<\varepsilon<\delta\rho_1-1$, there is a sufficiently small constant $0<c_1<1/2$ such that for all $0<y\le 2c_1$
	\[y^{\delta-1}(x+y^\delta)^{\rho_1-1}
	\le y^{\delta-1+\delta(\rho_1-1)}=y^{\delta\rho_1-1}
	\le (\rho_1\delta)^{-1}\]
	and
	\[
	(x+y^\delta)^{\rho_1-1}y^{\delta-2}\le
	y^{\delta\rho_1-2}
	=y^{\varepsilon-1}y^{\delta\rho_1-\varepsilon-1}
	\le
	4^{-1}\varepsilon[\rho_1(\delta-1)\delta]^{-1}y^{\varepsilon-1}.
	\]
	Thus, by \eqref{3.13b},
	\beqlb\label{3.18}
	h_y'(x,y)\ge-2\rho \hat{h}(x,y)^{\rho-1},~~
	h_{yy}''(x,y)
	\ge
	2^{-2}\rho\varepsilon \hat{h}(x,y)^{\rho-1}y^{\varepsilon-1},\qquad
	0<y\le 2(c_0\wedge c_1),
	\eeqlb
	where the constant $c_0$ is determined in \eqref{3.13b}.

	Recall the function $K(v,z)$ defined in \eqref{3.7bb}.
	By \eqref{3.8a}, for $0\le v\le1$ we have
	\beqnn
	K''_{zz}(v,z)
	\ge
	\rho\rho_1(1-\rho_1)v(1+z)^{\rho_1-2}\big[v[(1+z)^{\rho_1}-1]+1\big]^{\rho-1}
	\ge
	\rho\rho_1(1-\rho_1)v(1+z)^{\rho\rho_1-2}
	\eeqnn
	and then by \eqref{6.0a} and \eqref{3.0a},
	\beqlb\label{3.20aa}
	\ar\ar
	\int_0^\infty K(v,z)\mu_1(\dd z)
	=\int_0^\infty z^2\mu_1(\dd z) \int_0^1
	K''_{zz}(v,zu)(1-u)\dd u \cr
	\ar\ge\ar
	\rho\rho_1(1-\rho_1)v\int_0^\infty z^2\mu_1(\dd z) \int_0^1
	(1+zu)^{\rho\rho_1-2}(1-u)\dd u=
	\rho\rho_1(1-\rho_1)c_1(\rho\rho_1)v.
	\eeqlb
Now, by the arguments for the first equation
 of \eqref{3.14},
	\beqlb\label{3.20a}
	\int_0^\infty K_z^1h(x,y)\mu_1(\dd z)
	\ar=\ar
	-\hat{h}(x,y)^\rho (x+y^\delta)^{-\alpha_1}
	\int_0^\infty K(\frac{(x+y^\delta)^{\rho_1}}{\hat{h}(x,y)},z)\mu_1(\dd z) \cr
	\ar\ge\ar
	\rho\rho_1(1-\rho_1)c_1(\rho\rho_1)
	\hat{h}(x,y)^{\rho-1} (x+y^\delta)^{\rho_1-\alpha_1}.
	\eeqlb

	By \eqref{6.0a} and \eqref{3.18},
	for $0<y\le c_0\wedge c_1$
	\beqlb\label{3.19}
	\ar\ar
	\int_0^{c_0} K_z^2h(x,y)\mu_2(\dd z)
	=\int_0^{c_0} z^2\mu_2(\dd z) \int_0^1
	h_{yy}''(x,y+uz)(1-u)\dd u \cr
	\ar\ar\qquad\ge
	2^{-2}\rho\varepsilon
	\int_0^{c_0} z^2\mu_2(\dd z) \int_0^1\hat{h}(x,y+uz)^{\rho-1} (y+uz)^{\varepsilon-1}(1-u)\dd u.
	\eeqlb
	By \eqref{3.29},
for all $0<y\le c_0$ and $0<z,u\le1$,
	\beqnn
	\hat{h}(x,y+yuz)
= [x+(y+yuz)^\delta]^{\rho_1}+\tilde{h}(y+yuz)
	\le
	(x+2^\delta y^\delta)^{\rho_1}+4\tilde{h}(y)
	\le
	2^{\delta\rho_1+1}\hat{h}(x,y).
	\eeqnn
	Now replacing the variable $z$ with $yv$ in \eqref{3.19}
and using the fact $\mu_2(\dd (yv))=y^{-\alpha_2}\mu_2(\dd v)$, for $0<y\le c_0\wedge c_1$ we have
	\beqnn
	\ar\ar
	\int_0^{c_0} K_z^2h(x,y)\mu_2(\dd z) \cr
	\ar\ar\quad\ge
	2^{-2}\varepsilon\rho
	y^{2-\alpha_2}\int_0^{c_0/y} v^2\mu_2(\dd v) \int_0^1\hat{h}(x,y+yuv)^{\rho-1}(y+yuv)^{\varepsilon-1}
	(1-u)\dd u  \cr
	\ar\ar\quad\ge
	2^{\varepsilon-3}\varepsilon \rho
	y^{\varepsilon+1-\alpha_2}\int_0^1 z^2\mu_2(\dd z) \int_0^1\hat{h}(x,y+yuz)^{\rho-1}(1-u)\dd u \cr
	\ar\ar\quad\ge
	2^{\varepsilon-3+(\delta\rho_1+1)(\rho-1)}\varepsilon \rho
	y^{\varepsilon+1-\alpha_2} \hat{h}(x,y)^{\rho-1}\int_0^1 z^2\mu_2(\dd z)
	=:
	2^{-2}\rho \varepsilon\tilde{c}\hat{h}(x,y)^{\rho-1} y^{\varepsilon+1-\alpha_2}.
	\eeqnn
	On the other hand, for all $0<x,y\le c_0$,
	\beqnn
	\int_{c_0}^\infty K_z^2h(x,y)\mu_2(\dd z)
	\ge
	-\int_{c_0}^\infty \hat{h}(x,y+z)^\rho\mu_2(\dd z)
	\ge
	-\int_{c_0}^\infty \hat{h}(c_0,c_0+z)^\rho\mu_2(\dd z)
	=:-\tilde{c}_1
	\eeqnn
	with $\tilde{c}_1<\infty$ by the assumption of $\delta\rho\rho_1<1$.
	Then
	\beqlb\label{3.20b}
	\int_0^\infty K_z^2h(x,y)\mu_2(\dd z)
	\ge
	2^{-2}\rho \varepsilon\tilde{c} \hat{h}(x,y)^{\rho-1} y^{\varepsilon+1-\alpha_2}
	-\tilde{c}_1,\qquad 0<x,y\le c_0\wedge c_1.
	\eeqlb
	It thus follows from \eqref{3.0}, \eqref{3.20}--\eqref{3.18},
	\eqref{3.20a} and \eqref{3.20b} that
	\beqnn
	\ar\ar
	\mathcal{L}h(x,y) \cr
	\ar\ge\ar
	\rho \hat{h}(x,y)^{\rho-1}\Big[
	\rho_1(1-\rho_1)
	\big[b_{11}(x+y^\delta)^{\rho_1-2}x^{r_{11}}
	+ b_{21}c_1(\rho\rho_1) (x+y^\delta)^{\rho_1-\alpha_1}x^{r_{12}}\big] \cr
	\ar\ar\qquad \qquad\qquad
	+2^{-2} \varepsilon [b_{21}y^{r_{21}+\varepsilon-1}
	+b_{22}\tilde{c}y^{r_{22}+\varepsilon+1-\alpha_2}] \cr
	\ar\ar\qquad \qquad\qquad
	-a_1\rho_1x^{\theta_1}y^{\kappa_1}(x+y^\delta)^{\rho_1-1}
	-2a_2y^{\theta_2}x^{\kappa_2}\Big]
	-b_{22}\tilde{c}_1y^{r_{22}} \cr
	\ar\ge\ar
	\rho \hat{h}(x,y)^{\rho-1}[I_1(x,y)-I_2(x,y)]
	-b_{22}\tilde{c}_1y^{r_{22}} ,\qquad 0<x,y\le c_0\wedge c_1
	\eeqnn
	with $\tilde{c}_2:=\rho_1(1-\rho_1)[b_{11}\wedge({\blue b_{12}}c_1(\rho\rho_1))]$
	and
	$\tilde{c}_3:=2^{-2} \varepsilon [b_{21}\wedge(b_{22}\tilde{c})]$.

Now we are ready to conclude the proof of Lemma \ref{t7.5} by Lemma \ref{t2.3}(ii) again.

\noindent{\it Proof of Lemma \ref{t7.5}.}	
Since the proofs are similar, we only state that under \eqref{3.122a}.
	Under \eqref{3.122a},
	there is a constant $0<\rho_1<1$ satisfying
	\beqlb\label{4.1aa}
	\rho_1<(1-\theta_1)\wedge(2-\kappa_2)
	\eeqlb
	and
	\beqnn
	\rho_1>\frac{1-\theta_1}{\kappa_1-r_2},\quad
	\rho_1<\frac{\kappa_2}{r_2+1-\theta_2},\quad
	\rho_1<\frac{\kappa_2-r_1}{1-\theta_2},
	\eeqnn
	which implies that there is a constant $\delta>\rho_1^{-1}>1$ such that
	\beqnn
	\delta(1-\rho_1-\theta_1)<\kappa_1-1-r_2,\quad
	r_2+1-\theta_2-\delta\kappa_2<0
	\eeqnn
	and
	\beqlb\label{4.1bb}
	r_1<\kappa_2-\rho_1+\theta_2/\delta,
	\eeqlb
	respectively.
	Moreover, there is a small enough constant $0<\varepsilon<\delta\rho_1-1$ such that
	$r_2+\varepsilon<0$,
	\beqlb\label{4.1a}
	\delta(1-\rho_1-\theta_1)<\kappa_1-1-r_2-\varepsilon
	\eeqlb
	and
	\beqlb\label{4.1b}
	r_2+\varepsilon+1-\theta_2-\delta\kappa_2<0.
	\eeqlb
	Let $0<\rho<1$ be small enough satisfying
	\beqlb\label{4.1c}
	\rho\rho_1+r_1<0,\qquad r_2+\rho+\varepsilon<0
	\eeqlb
	and $\delta\rho\rho_1<1$.
	
	Let $I_1$ and $I_2$ be the functions defined in Lemma \ref{t7.5b}.
	As $x^{\rho_1}\le y\le 1$, by Lemma \ref{t2.3}(ii) and the fact
	$\tilde{h}(y)\le y$ for $y>0$ we get
	\beqnn
	\hat{h}(x,y)\le x^{\rho_1}+y^{\delta\rho_1}+y
	\le 3y
	\eeqnn
	and then
	\beqnn
	\hat{h}(x,y)^{\rho-1}I_1(x,y)
	 \ar\ge\ar
	\tilde{c}_3\hat{h}(x,y)^{\rho-1} y^{\varepsilon+1}[y^{r_{21}-2}+y^{r_{22}-\alpha_2}] \cr
	 \ar\ge\ar
	3^{\rho-1}\tilde{c}_3y^{\varepsilon+\rho}[y^{r_{21}-2}+y^{r_{22}-\alpha_2}]
	=3^{\rho-1}\tilde{c}_3
	y^{\varepsilon+\rho}[y^{r_{21}-2}+y^{r_{22}-\alpha_2}].
	\eeqnn
	Since $y\le x^{\rho_1}\le1$, we get
	\beqnn
	\hat{h}(x,y)\le x^{\rho_1}+y^{\delta\rho_1}+y
	\le x^{\rho_1}+2y\le3x^{\rho_1},~~
	x+y^\delta\le x+x^{\rho_1\delta}\le 2x
	\eeqnn
	and then
	\beqnn
	\hat{h}(x,y)^{\rho-1}I_1(x,y)
	\ar\ge\ar
	\tilde{c}_2\hat{h}(x,y)^{\rho-1}\big[(x+y^\delta)^{\rho_1-2}x^{r_{11}}
	+(x+y^\delta)^{\rho_1-\alpha_1}x^{r_{12}}\big] \cr
	\ar\ge\ar
	3^{\rho-1}\tilde{c}_2x^{\rho\rho_1-\rho_1}
	\big[2^{\rho_1-2}x^{\rho_1-2+r_{11}}
	+2^{\rho_1-\alpha_1}x^{\rho_1-\alpha_1+r_{12}}\big].
	\eeqnn
	It follows from \eqref{1.6} that
	\beqlb\label{a3.5}
	\hat{h}(x,y)^{\rho-1}I_1(x,y)\ge
	3^{\rho-1}2^{\rho_1-2}(\tilde{c}_2\wedge\tilde{c}_3)
	[x^{\rho\rho_1+r_1}\wedge y^{r_2+\rho+\varepsilon}],\qquad 0<x,y\le1.
	\eeqlb

	By \eqref{4.1aa} and \eqref{4.1a}, there is a constant $0<c_2\le c_0\wedge c_1$
	such that
	\beqlb\label{a3.1}
	\tilde{c}_3y^{r_2+\varepsilon+1}
	\ar=\ar
	\tilde{c}_3 y^{r_2+\varepsilon+1-\delta(\rho_1-1+\theta_1)}y^{\delta(\rho_1-1+\theta_1)}
	\ge
	\tilde{c}_3 y^{r_2+\varepsilon+1-\delta(\rho_1-1+\theta_1)-\kappa_1}(x+y^{\delta})^{\rho_1-1+\theta_1}
	y^{\kappa_1} \cr
	\ar\ge\ar
	2a_1\rho_1 x^{\theta_1} y^{\kappa_1}(x+y^\delta)^{\rho_1-1},\qquad 0<y\le c_2.
	\eeqlb
	Similarly, by \eqref{4.1aa}, \eqref{4.1bb} and \eqref{4.1b}, there is a constant $0<c_3\le c_2$
	such that for all $0<x,y\le c_3$
	\beqnn
	\ar\ar
	\tilde{c}_3 y^{r_2+\varepsilon+1}(x+y^\delta)^{2-\rho_1}
	\ge
	\tilde{c}_3  y^{r_2+\varepsilon+1-\theta_2} y^{\theta_2}x^{\kappa_2} (x+y^{\delta})^{2-\rho_1-\kappa_2} \cr
	\ar\ar\quad\ge
	\tilde{c}_3 y^{r_2+\varepsilon+1-\theta_2} y^{\theta_2}x^{\kappa_2} y^{\delta(2-\rho_1-\kappa_2)} =
	\tilde{c}_3  y^{r_2+\varepsilon+1-\theta_2-\delta\kappa_2} y^{\theta_2}x^{\kappa_2} y^{\delta(2-\rho_1)} \cr
	\ar\ar\quad\ge
	4a_2 2^{1-\rho_1+\theta_2/\delta} x^{\kappa_2}y^{\delta(2-\rho_1+\theta_2/\delta)}
	\eeqnn
	and
	\beqnn
	\tilde{c}_2[x^{r_{11}}+(x+y^\delta)^{2-\alpha_1}x^{r_{12}}]
	\ar\ge\ar
	\tilde{c}_2x^{r_1+2}
	=
	\tilde{c}_2x^{r_1+2-(\kappa_2+2-\rho_1+\theta_2/\delta)} x^{\kappa_2} x^{2-\rho_1+\theta_2/\delta} \cr
	\ar\ge\ar
	4a_22^{1-\rho_1+\theta_2/\delta}x^{\kappa_2} x^{2-\rho_1+\theta_2/\delta}.
	\eeqnn
	Thus, by Lemma \ref{t2.3}(ii) again,
	\beqlb\label{a3.2}
	\ar\ar
	I_1(x,y)(x+y^\delta)^{2-\rho_1}
	\ge
	4a_2 2^{1-\rho_1+\theta_2/\delta} x^{\kappa_2} x^{2-\rho_1+\theta_2/\delta}
	+
	4a_2 2^{1-\rho_1+\theta_2/\delta}x^{\kappa_2} y^{\delta(2-\rho_1+\theta_2/\delta)} \cr
	\ar\ar\quad\ge
	4a_2 x^{\kappa_2} (x+y^\delta)^{2-\rho_1+\theta_2/\delta}
	\ge
	4a_2 y^{\theta_2}x^{\kappa_2} (x+y^\delta)^{2-\rho_1},\quad
	0<x,y\le c_3.
	\eeqlb
	Combining \eqref{a3.1} and \eqref{a3.2} we get
	$2^{-1}I_1(x,y)\ge I_2(x,y)$ for all $0<x,y\le c_3$,
	which together with Lemma \ref{t7.5b} and
	\eqref{a3.5} implies
	\beqnn
	\mathcal{L}h(x,y)
	\ge
	3^{\rho-1}2^{\rho_1-3}(\tilde{c}_2\wedge\tilde{c}_3)
	[x^{\rho\rho_1+r_1}\wedge y^{r_2+\rho+\varepsilon}]
	-b_{22}\tilde{c}_1y^{r_{22}},\qquad 0<x,y\le c_3.
	\eeqnn
	By \eqref{4.1c}, there are constants $C>0$
	and $0<c\le c_3$ such that $\mathcal{L}h(x,y)\ge C$ for all
	$0<x,y\le c$.
	\qed

In the following we introduce the key test
 function $h$ and the estimate on $\mathcal{L}h$ under the assumptions of Theorem \ref{t1.2}(iii).
	\blemma\label{t7.7}
Recall the function $\tilde{h}$ defined by \eqref{7.5}
	for $0<\varepsilon<2^{-1}\wedge(r_2-1-\kappa_1)\wedge(r_1-1-\kappa_2)$.
	For $\theta_1,\theta_2>0$ let $\hat{h}_{\theta_1,\theta_2}(x,y):=x^{1-\theta_1}+y^{1-\theta_2}$;
	for $\theta_1>0$ and $\theta_2=0$ let $\hat{h}_{\theta_1,0}(x,y):=x^{1-\theta_1}+\tilde{h}(y)$;
	for $\theta_1=0$ and $\theta_2>0$ let $\hat{h}_{0,\theta_2}(x,y):=y^{1-\theta_2}+\tilde{h}(x)$;
	for $\theta_1=\theta_2=0$
	let $\hat{h}_{0,0}(x,y):=\tilde{h}(x)+\tilde{h}(y)$.
	For $0<\rho<1$ let $h(x,y):=- \hat{h}_{\theta_1,\theta_2}(x,y)^\rho$.
	Then, under the assumptions of Theorem \ref{t1.2}(iii),
	there are constants $C,c,\varepsilon>0$ such that
	$\mathcal{L}h(x,y)\ge C$ for all $0<x,y\le c$.
	\elemma
	\proof
The proof is similar to those of Lemmas \ref{t7.2} and \ref{t7.4}.  We leave the details to interested readers.
	\qed
	
	Now we are ready to complete the proofs of Theorems \ref{t1.2}
and \ref{t1.4bb}.
	
	\noindent{\it Proof of Theorem \ref{t1.2}.}
	Under the assumptions in Theorem \ref{t1.2}(i),
	there are constants $\rho,\rho_1,\rho_2$ determined
	in Lemma \ref{t7.2} and
	let the function $h$ be defined by \eqref{3.7}.
	Under the assumptions of Theorem \ref{t1.2}(ii),
	let $h$ be the function defined in Lemmas \ref{t7.4} and \ref{t7.5}.
	For the assumptions of Theorem \ref{t1.2}(iii),  let the function $h$ be determined in Lemma \ref{t7.7}.
	Let the constants $C,c>0$ be the constants determined in Lemmas \ref{t7.2}, \ref{t7.4}, \ref{t7.5} and \ref{t7.7}.
	Set 
	\beqlb\label{7.2}
	g(x,y):=|h(c,0)|\wedge |h(0,c)|+h(x,y),
	\qquad x,y>0.
	\eeqlb
	Then the condition (i) of Proposition \ref{t2.2} is obvious
	and $g(x,y)\le 0$ for all $x\ge c$ or $y\ge c$.
	By Lemma \ref{t7.2}, \ref{t7.4}, \ref{t7.5} and \ref{t7.7}, for all $0<x,y\le c$,
	\beqnn
	\mathcal{L}g(x,y)\ge C
	\ge C[|h(c,0)|\wedge |h(0,c)|]^{-1} g(x,y).
	\eeqnn
	Then by Proposition \ref{t2.2}, one gets $\mbf{P}\{\tau_0<\infty\}\ge
g(X_0,Y_0)[|h(c,0)|\wedge |h(0,c)|]^{-1}$ for all small enough
	$X_0,Y_0>0$.
	For general, $X_0,Y_0>0$, applying the strong Markov property and the same arguments
	at the end of the proof of \cite[Theorem 1.3]{XYZh24} we conclude the proof.
	\qed

	\noindent{\it Proof of Theorem \ref{t1.4bb}.}
	Under the assumptions,
	let  $\rho,\rho_1,\rho_2,\delta_0$ be the constants given
	in Lemma \ref{t7.3} and
	 $h:=\tilde{h}$ be the function given by \eqref{3.7b}.
By the same arguments as in the proof of Theorem \ref{t1.2},
	one completes the proof.
	\qed
	
\bremark
By the proofs of Theorems \ref{t1.2a}, \ref{t1.2} and  \ref{t1.4bb},
there are constants $C,c>0$ such that
$\mbf{P}\{\tau_0<\infty\}\ge C$ for all $0<X_0,Y_0\le c$ under the assumptions
in Theorems \ref{t1.2a}, \ref{t1.2} and  \ref{t1.4bb}.
\eremark

\end{document}